\documentclass[10pt]{amsart}

\usepackage{amsmath}
\usepackage{amsthm}
\usepackage{amssymb}
\usepackage{amscd}
\usepackage[plainpages=false,colorlinks,hyperindex,pdfpagemode=None,bookmarksopen,linkcolor=red,citecolor=blue,urlcolor=blue]{hyperref}

\swapnumbers
\theoremstyle{plain}
\newtheorem{theorem}[subsubsection]{Theorem}
\newtheorem*{theorem*}{Theorem}
\newtheorem{proposition}[subsubsection]{Proposition}
\newtheorem{lemma}[subsubsection]{Lemma}
\newtheorem{corollary}[subsubsection]{Corollary}
\newtheorem*{conjecture}{Conjecture}
\theoremstyle{definition}
\newtheorem*{definition}{Definition}
\theoremstyle{remark}
\newtheorem*{remark}{Remark}
\newtheorem{example}[subsubsection]{Example}
\newtheorem*{acknowledgements}{Acknowledgements}

\renewcommand{\AA}{\mathbf{A}}
\newcommand{\BB}{\mathbf{B}}
\newcommand{\CC}{\mathbb{C}}
\newcommand{\DD}{\mathbf{D}}
\newcommand{\EE}{\mathbf{E}}
\newcommand{\FF}{\mathbf{F}}
\newcommand{\GG}{\mathbf{G}}
\newcommand{\HH}{\mathbf{H}}
\newcommand{\LL}{\mathbf{L}}
\newcommand{\MM}{\mathbf{M}}

\newcommand{\PP}{\mathbf{P}}
\renewcommand{\SS}{\mathbf{S}}
\newcommand{\TT}{\mathbf{T}}
\newcommand{\UU}{\mathbf{U}}
\newcommand{\VV}{\mathbf{V}}
\newcommand{\XX}{\mathbf{X}}
\newcommand{\YY}{\mathbf{Y}}
\newcommand{\ZZ}{\mathbf{Z}}

\newcommand{\ur}{\textrm{ur}}

\newcommand{\Ind}{\rm{Ind}}
\newcommand{\Hom}{\rm{Hom}}
\newcommand{\End}{\rm{End}}
\newcommand{\Aut}{\rm{Aut}}
\newcommand{\varchi}{\mathcal{X}}
\newcommand{\Gm}{\mathbb{G}_m}
\newcommand{\Ga}{\mathbb{G}_a}
\newcommand{\supp}{\rm{supp}}
\newcommand{\GGL}{\rm{\mathbf{GL}}}
\newcommand{\GL}{\rm{GL}}
\newcommand{\PPGL}{\rm\mathbf{\mathbf{PGL}}}
\newcommand{\PGL}{\rm{PGL}}
\newcommand{\SSL}{\rm{\mathbf{SL}}}
\newcommand{\SL}{\rm{SL}}
\newcommand{\Ssp}{\rm{\mathbf{Sp}}}
\newcommand{\geom}{\rm{geom}}
\newcommand{\Gal}{\rm{Gal}}

\newcommand{\tr}{\rm{tr}}
\newcommand{\Adm}{\mathit{Adm}}
\newcommand{\rk}{\operatorname{rk}}
\newcommand{\codim}{\operatorname{codim}}

\begin{document}
\setcounter{tocdepth}{2}
\title[On the unramified spectrum of spherical varieties]{On the unramified spectrum of spherical varieties over $p$-adic fields} 
\author{Yiannis Sakellaridis} \email{yiannis@post.tau.ac.il}
\thanks{Partially supported by the Marie Curie Research Training Network in Arithmetic Algebraic Geometry, Contract MRTN-CT-2003-504917.}

\address{School of Mathematical Sciences, Tel Aviv University, Ramat Aviv, Tel Aviv 69978, Israel.}

\begin{abstract}
The description of irreducible representations of a group $G$ can be seen as a question in harmonic analysis; namely, decomposing a suitable space of functions on $G$ into irreducibles for the action of $G \times G$ by left and right multiplication.

For a split $p$-adic reductive group $G$ over a local non-archimedean field, unramified irreducible smooth representations are in bijection with semisimple conjugacy classes in the ``Langlands dual'' group.

We generalize this description to an arbitrary spherical variety $X$ of $G$ as follows: Irreducible unramified quotients of the space $C_c^\infty(X)$ are in natural ``almost bijection'' with a number of copies of $A_X^*/W_X$, the quotient of a complex torus by the ``little Weyl group'' of $X$. This leads to a description of the Hecke module of unramified vectors (a weak analog of geometric results of Gaitsgory and Nadler), and an understanding of the phenomenon that representations ``distinguished'' by certain subgroups are functorial lifts. In the course of the proof, rationality properties of spherical varieties are examined and a new interpretation is given for the action, defined by F.~Knop, of the Weyl group on the set of Borel orbits.
\end{abstract}



\maketitle

\tableofcontents

\section{Introduction}\label{secintro}

\subsection{Motivation} Let $\GG$ be a split reductive linear algebraic group over a local non-archimedean field $k$ in characteristic zero. 
A $k$-variety $\XX$ with a $k$-action\footnote{Our convention will be that the action of the group is on the right.} of $\GG$ is called \emph{spherical} if the Borel subgroup $\BB$ of $\GG$ has an open orbit $\mathring \XX$ on $\XX$.\footnote{It is enough to assume that there exists an open Borel orbit over the algebraic closure; we show later that it will then have a point over $k$.} This includes, but is not limited to, symmetric and horospherical varieties (these are homogeneous spaces whose isotropy groups are, respectively, equal to the subgroup of points fixed by an involution of $\GG$, or containing a maximal unipotent subgroup $\UU$). The group $\GG$ itself can be considered as a spherical variety under the action of $\GG\times \GG$ on the left and right, and in fact many well-known theorems for algebraic groups can be seen as special cases of more general theorems for spherical varieties under this perspective (e.g.\ \cite{KnHC}).\footnote{\label{footnote1}Notice, though, that one usually makes use of a theorem for $\GG$ in order to prove its generalization to an arbitrary spherical variety; this is the case in our present work, too.}

The importance of the open orbit condition becomes apparent in the following:
\begin{theorem*}[Vinberg and Kimel'feld \cite{VK}]
Let $\XX$ be a quasi-affine $\GG$-variety over an algebraically closed field $k$. The space $k[\XX]$ of regular functions on $\XX$, considered as a representation of $\GG$ by right translations, is multiplicity-free if and only if $\XX$ is spherical.
\end{theorem*}
Hence, if $\XX=\HH\backslash\GG$ is quasi-affine (which, as we discuss in \S \ref{ssspherical}, can be assumed without serious loss of generality) the above theorem states that $\XX$ is spherical if and only if $(\GG,\HH)$ is a \emph{Gelfand pair} in the category of algebraic representations.

One goal of the present work is to examine to what extent a similar result is true in the category of smooth representations of $p$-adic groups, and more generally to describe (part of) the representations contained in the spectrum of a spherical variety over a $p$-adic field. More precisely, we consider the unramified component of $C_c^\infty(X)$, the representation of $G=\GG(k)$ on the space of smooth, compactly supported functions on  $X=(\HH\backslash\GG)(k)$, and provide a description of the (generic) irreducible quotients of this representation. Equivalently, this amounts to a description of embeddings $\pi\hookrightarrow C^\infty(X)$, for $\pi$ an irreducible unramified representation in general position, or in the simple case where $X=H\backslash G$, to a description of the space of $H$-invariant functionals on $\pi$. If $H$ contains the unipotent radical $U_P$ of a parabolic, and $\Psi$ is an unramified character of $U_P$ normalized by $H$, then we also consider sections of the corresponding line bundle $\mathcal L_\Psi$ over $H\backslash G$ (for instance, the Whittaker model); for simplicity we will mostly ignore this case in the introduction and refer the reader to \ref{ssparindthm}.

The description itself leads to a fascinating picture, in which a ``dual group'' $\widehat{G}_X$, a reductive subgroup of the ``Langlands dual'' group of $G$, seems to play a role in parametrizing the irreducible representations appearing in the spectrum of $X$. This dual group is implicit in the purely algebro-geometric work of Knop \cite{KnAu}, it has appeared in recent work of Gaitsgory and Nadler \cite{GN1,GN2} in the context of the Geometric Langlands Program, it appears in our present work on the unramified spectrum and conforms to the philosophy established by the deep work of several people on particular cases, in both the local and automorphic settings. (More details below.) It is natural to ask to what extent and in precisely what fashion it plays a role in describing the whole spectrum of $X$.\footnote{In fact, the global Langlands conjecture with a reasonable version of an ``automorphic Langlands group'', together with the Burger-Sarnak principle, our results and a Chebotarev density argument could ``in most cases'' imply such a role, but since all this is anyway speculative we will not discuss this issue further.}

Spherical varieties are ubiquitous in the theory of automorphic forms, although their applications there have never been examined in this generality. Providing candidates for Gelfand pairs (but even, sometimes, when the Gelfand condition fails), spherical varieties play an essential role in the theory of integral representations of $L$-functions \cite{PS, GPR}, in the relative trace formula \cite{Ja, La} and other areas such as explicit computations of arithmetic interest \cite{Csph, HiHe, Of}. 

On the other hand, the theory of spherical varieties has been greatly developed from the algebro-geometric point of view in the work of Brion, Knop, Luna, Vinberg, Vust and others. They have discovered rich geometric and combinatorial structures related to spherical varieties, and these structures will lend us the dictionary for describing the spectrum. Hence, another goal of the present work is to initiate a systematic study of the representation theory of spherical varieties by establishing a connection with these algebraic structures. This allows to replace explicit, hands-on methods such as double coset decompositions with more elegant ones and leads to the general picture to which I alluded above. In the course of establishing this connection, I have had to examine several rationality properties of the varieties over $k$, which may be of independent interest. 

For the rest of the introduction I describe more precisely some of the results and comment on the methods of proof.

\subsection{Parametrization of irreducible quotients} The main phenomenon that the current work reveals is the local analog of a global statement of the following form, very often arising in the theory of the relative trace formula and elsewhere: ``An automorphic representation $\pi$ of $\GG$ is a functorial lift from (a certain other group) $\GG'$ if and only if it is \emph{distinguished} by (a certain subgroup) $\HH$.'' Instead of explaining the global notion of being ``distinguished'', we describe its local analog which is the object of study here: $\pi$ is \emph{distinguished} by $H=\HH(k)$ if it appears as a quotient of $C_c^\infty(X)$, the representation of $G=\GG(k)$ on the space of smooth, compactly supported functions on  $X=(\HH\backslash\GG)(k)$. 

Recall that irreducible unramified representations of $G$ are in ``almost bijection'' with semisimple conjugacy classes ($\simeq A^*/W$) in the Langlands dual group $\widehat{G}$ of $G$ -- and they can be realized as (subquotients of) unramified principal series $I(\chi)=\Ind_B^G(\chi\delta^{\frac{1}{2}})$. (Here $A^*$ is the maximal torus of the dual group, or equivalently the complex torus of unramified characters of the Borel subgroup, and $W$ is the Weyl group; $\delta$ is the modular character of the Borel subgroup.) To each spherical variety $\XX$, Brion \cite{BrGe} associates a finite group (the ``little Weyl group'') $W_X$ acting faithfully on the vector space $\mathfrak a_{X,\mathbb Q}^*:= \varchi (\mathring\XX) \otimes_{\mathbb Z} \mathbb Q$, where $\varchi(\mathring\XX)$ denotes the weights of $\BB$-semiinvariants (regular eigenfunctions) on the open orbit. Let us denote by $\PP(\XX)$ the standard parabolic $\{g\in \GG| \mathring \XX g = \mathring \XX\}$ and by $[W/W_{P(X)}]$ the canonical set of representatives of minimal length of $W/W_{P(X)}$-cosets (where $W_{P(X)}$ is the Weyl group of the Levi of $\PP(\XX)$). An alternative construction by Knop \cite{KnWe,KnAs} proves, among others, that $W_X\subset [W/W_{P(X)}]$ canonically; and the action of $W_X$ on $a_{X,\mathbb Q}^*$ is generated by reflections. The complex analog of $\mathfrak a_{X,\mathbb Q}^*$ is the Lie algebra $\mathfrak a_X^*:= \varchi (\mathring\XX) \otimes_{\mathbb Z} \CC$ of a subtorus $A_X^*\subset A^*$.  We assume throughout (without serious loss of generality, cf.\ \S\ref{ssspherical}, \ref{invariantmeasure}) that $\XX$ is quasi-affine and $\mathring X$ (the set of $k$-points of its open Borel orbit) carries a $B$-invariant measure. A simplified version of our main result is (cf.\ Theorems \ref{existence}, \ref{mainthm}):

\begin{theorem}\label{Theorem1}
A necessary condition for the existence of a non-zero morphism: $C_c^\infty(X)\to I(\chi)$ is that: 
$$\chi\in {^w\!\left(\delta^{-\frac{1}{2}}A_X^*\right)} \textrm{ for some } w\in [W/W_{P(X)}].$$

If $\chi\in \delta^{-\frac{1}{2}}A_X^*$ then $\chi\delta^{\frac{1}{2}}$ is a character of $P(X)$ and almost all unramified irreducible $\pi$ admitting a nonzero morphism from $C_c^\infty(X)$ are isomorphic to $\Ind_{P(X)}^G (\chi\delta^\frac{1}{2})$ for such a $\chi$. 

Moreover, for almost every such $\pi$ we have:
$$\dim\Hom(C_c^\infty(X),\pi)=(\mathcal N_W(\delta^{-\frac{1}{2}} A_X^*):W_X)\times{|H^1(k,\AA_X)|}$$
where $\AA_X$ is the image in $\BB/\UU$ of the stabilizer of a generic point on $\XX$. (The factor $|H^1(k,\AA_X)|$ is equal to the number of $B$-orbits on $\mathring X$.)
\end{theorem}

The word ``almost every'' refers to the variety structure of $A^*/W$ and the subspaces under consideration, and means ``except for a subvariety of strictly smaller dimension''. 

In other words, we describe a natural ``almost one-to-one correspondence'' between a basis of irreducible quotients of the ``unramified'' Bernstein component (cf.\ \S \ref{ssBernstein}) $\left(C_c^\infty(X)\right)_\ur$ and $|H^1(k,\AA_X)|$ copies of the complex space $\delta^{-\frac{1}{2}}A_X^*/W_X$. In the phenomenon of ``distinguished'' lifts that I alluded to before, $A_X^*$ and $W_X$ are the maximal torus and the Weyl group of the ``Langlands dual'' group $\widehat{G}_X\subset \widehat G$, and distinguished unramified representations (and conjecturally, not only unramified) should be functorial lifts from some group $G'$ with $\widehat{G'}=\widehat G_X$.\footnote{Note, however, that global distinguishedness typically involves a condition additional to local functoriality from $G'$ at every place, which is usually expressed by the non-vanishing of an $L$-value. This is a more complicated problem that will not concern us here.}

\subsection{The Hecke module of unramified vectors} The dual group $\widehat G_X$ has appeared in a more canonical way in recent work of Gaitsgory and Nadler \cite{GN1,GN2} in the context of the geometric Langlands program.\footnote{The Weyl group of the dual group of Gaitsgory and Nadler has not been identified yet in their work, but it is conjectured to be equal to $W_X$.} In that work, it is proven that a certain category of $G(\mathfrak o)$-equivariant perverse sheaves on (a quasi-global analog of) the space of $k$-points on $\XX$, where the spherical variety $\XX$ is now defined over a global complex curve, $\mathfrak o=\CC[[t]]$ and $k=\CC((t))$, is equivalent to the category of finite-dimensional representations of $\widehat{G}_X$. I prove the following weak analog of their results in the $p$-adic setting:

\begin{theorem}\label{Theorem4}
Let $K$ denote a hyperspecial maximal compact subgroup of $G$. Let $\mathcal H_X$ denote the quotient of the Hecke algebra $\mathcal H(G,K)\simeq \CC[A^*]^W$ corresponding to the image of $\delta^{-\frac{1}{2}}A_X^*$ under $A^*\to A^*/W$ and let $\mathcal K_X$ denote the quotient field of $\mathcal H_X$.

The space $C_c^\infty(X)^K$ is a finitely-generated, torsion-free module for $\mathcal H_X$. 

Moreover, we have: $C_c^\infty(X)^K\otimes_{\mathcal H_X} \mathcal K_X \simeq \left(\CC(\delta^{-\frac{1}{2}}A_X^*)^{W_X}\right)^{|H^1(k,\AA_X)|}$.
\end{theorem}

Notice that the invariants $A_X^*, W_X, \AA_X$ appearing in the theorems above only depend on the open $\GG$-orbit $\HH\backslash\GG\subset\XX$, although the representations considered depend on $\XX$ itself. The basic fact leading to this conclusion is that the representation-theoretic content of smaller $\GG$-orbits can be read off from the open orbit, more precisely certain intertwining operators supported on the smaller orbits appear as residues of operators on the open orbit (cf.\ Proposition \ref{pickup}).

As is usually the case with spherical varieties, one recovers classical results by considering $\GG$ as a spherical $\GG\times\GG$ variety under left and right multiplication; in this case, we recover the (generic) description of irreducible unramified representations of $G$ by semisimple conjugacy classes in its Langlands dual group, and the Satake isomorphism (without reproving them, since they are used throughout).

\subsection{Interpretation of Knop's action}
Theorem \ref{Theorem1} comes as a corollary to an analysis that we perform and an interpretation in the context of the representation theory of $p$-adic groups that we give to an action, defined by F.~Knop \cite{KnOrbits}, of the Weyl group $W$ on the set of Borel orbits on $\XX$. (Moreover, the case of a non-trivial character $\Psi$ on the unipotent radical $U_P\subset H$ of a parabolic is treated with the help of Knop's extension of this action to the non-spherical case.) We recall the definition of this action in Section \ref{secac}; for now, let $^w \YY$ denote the image of a $\BB$-orbit $\YY$ under the action of $w\in W$.
Using standard ``Mackey theory'' \cite{BZ, Cas} we define for every $B$-orbit of maximal rank $Y$ (for the definition of rank, cf.\ \ref{ssspherical}) a rational family of morphisms: $S_\chi^Y:C_c^\infty(X)\to I(\chi)$, given by rational continuation of a suitable integral on $Y$; what was denoted by $S_\chi$ in the formulation of Theorem \ref{Theorem1} is now $S_\chi^{\mathring X}$. (Notice that, while this form of ``Mackey theory'' has been used extensively in the past, it has probably never been applied in this generality, and the technical results that we collect or prove for that purpose may be of independent interest.) For simplicity, let us assume here that $H^1(k,\AA_X)=1$.

Recall now that we also have the standard intertwining operators $T_w: I(\chi)\to I({^w\chi})$.
The heart of the current work is the proof of the following theorem on the effect of composing the operators $T_w$ with $S_\chi^Y$:

\begin{theorem}
\label{Theorem2} 
We have $T_w \circ S_\chi^{\mathring X}\ne 0$ if and only if $w\in [W/W_{P(X)}]$. In that case, $T_w \circ S_\chi^{\mathring X} \sim S_{^w\chi}^{^w\!\mathring X}$, where $\sim$ denotes equality up to a non-zero rational function of $\chi$.
\end{theorem}

In the case where $H^1(k,\AA_X)\ne 1$, there are several morphisms associated to each $\BB$-orbit, and the theorem still holds for the spaces spanned by these morphisms. Moreover, for $w=w_\alpha$, a simple reflection, we exhibit explicit bases for these spaces which are mapped to each other by composition with $T_w$ (cf.\ Theorem \ref{maintool}).

\subsection{Rationality results}
There are several rationality properties of our varieties that need to be established in order to apply the algebro-geometric theory (developed over an algebraically closed field) to our problems. The key feature here is that since the group is split and the variety spherical, every $\BB$-eigenfunction on $\XX$ is (up to a constant) defined over $k$ (Lemma \ref{rationalorbit}). It is expected that in the non-split case the situation will be significantly more complicated. I summarize some of the rationality results proven in section \ref{secrat}:

\begin{theorem}
\label{Theorem3} 
\begin{enumerate}
 \item Every $\BB$-orbit of maximal rank has a point over $k$.
 \item If $\YY$ is a $\BB$-orbit defined over $k$, then so is $^w\YY$, for every $w\in W$.
 \item Every $\GG$-orbit is defined and has a point over $k$.
\end{enumerate}
\end{theorem}

The rationality properties of the structure of Borel orbits on $X$ that we examine may be of independent interest. They generalize a portion of work of Helminck and Wang on symmetric varieties \cite{HW}.

\subsection{Endomorphisms} Finally, I discuss a (partly conjectural) ring of endomorphisms of the Hecke module $(C_c^\infty(X))^K$ which bears a remarkable similarity to the algebra of invariant differential operators on a spherical variety. For simplicity, assume here that $H^1(k,\AA_X)=1$. We recall Knop's generalization of the Harish-Chandra homomorphism \cite{KnHC}:

\begin{theorem*}[Knop]
The algebra of invariant differential operators on a spherical variety $\XX$ (over an algebraically closed field $k$ of characteristic 0) is commutative and isomorphic to $k[\rho+\mathfrak a_X^*]^{W_X}$.
This generalizes the Harish-Chandra homomorphism for the center $\mathfrak z(G)$ of the universal enveloping algebra of $\mathfrak g$ (if we regard the group $\GG$ as a spherical $\GG\times \GG$ variety) and the following diagram is commutative:
\begin{equation*} \begin{CD}
\mathfrak{z}(\GG) @>>> \mathfrak{D}(\XX)^\GG \\
@|                       @| \\
k[\mathfrak a^*]^W @>>> k[\rho+\mathfrak a_X^*]^{W_X}
\end{CD}
\end{equation*}
\end{theorem*}

Our description of unramified vectors in $C_c^\infty(X)$ leads easily to a conjectural description of a commutative subalgebra of their endomorphism algebra as an $\mathcal H(G,K)$-module, which should be naturally isomorphic the $\CC[\delta^{-\frac{1}{2}}A_X^*]^{W_X}$, the ring of regular functions on $\delta^{-\frac{1}{2}}A_X^*/W_X$. The precise statement of the conjecture is:

\begin{conjecture}Call ``geometric'' an endomorphism of $(C_c^\infty(X))^K$ that preserves up to a rational multiple the family of morphisms $S_\chi^{\mathring X}$.
There is a canonical isomorphism $\left(\End_{\mathcal H(G,K)}C_c^\infty(X)^K\right)^{\geom} \simeq \CC[\delta^{-\frac{1}{2}}A_X^*]^{W_X}$ such that the following diagram commutes:
\begin{equation*} \begin{CD}
\mathcal H(G,K) @>>> \left(\End_{\mathcal H(G,K)}C_c^\infty(X)^K\right)^{\geom} \\
@|                @|  \\
\CC[A^*]^W @>>> \CC[\delta^{-\frac{1}{2}}A_X^*]^{W_X}
\end{CD}
\end{equation*}
\end{conjecture}

In fact, it is easy to prove this conjecture in many cases:

\begin{theorem}
The above conjecture is true if: 
\begin{enumerate}
\item
the unramified spectrum of $X$ is generically multiplicity-free, in which case the geometric endomorphisms are all the endomorphisms of $(C_c^\infty(X))^K$, or
\item
the spherical variety $X$ is ``parabolically induced'' from a spherical variety whose unramified spectrum is generically multiplicity-free.
\end{enumerate}
\end{theorem}

\subsection{Some notation} We will be working throughout over fields of characteristic zero. Unless otherwise stated, we denote by $k$ a local non-archimedean field (by which we mean a locally compact one -- hence, with finite residue field of order $q$), by $\mathfrak o$ its ring of integers, by $\varpi$ a uniformizing element and by $|\bullet|$ the standard $p$-adic absolute value on $k$. We generally denote by $\GG$ a \emph{reductive group}, that is an affine algebraic group with trivial unipotent radical; in addition, for the whole paper ``reductive'' will also mean (geometrically) connected, and the group will be split over the field of definition.

Given a scheme $\YY$ over $k$, we denote by $\YY_{\bar k}$ the base change $\YY\times_{\rm{spec}\,k} \rm{spec}\,\bar k$. The set of $k$-points will be denoted by $Y$ or by $\YY(k)$. 

We generally fix a Borel subgroup $\BB\subset \GG$ (with unipotent radical $\UU$) and a maximal torus $\AA\subset\BB$ (which, of course will also be identified with the reductive quotient of $\BB$). We also fix the corresponding root system and choice of positive roots. We denote by $\Gm,\, \Ga$ the multiplicative and additive groups, respectively, over $k$, by $\mathcal N(\bullet)$ the normalizer of $\bullet$, by $\mathcal L(\bullet)$ the Lie algebra of $\bullet$ and by $\UU_\bullet$ the unipotent radical of a group $\bullet$. If $\alpha$ is a root of $\AA$, $\UU_\alpha$ will denote the corresponding one-parameter unipotent subgroup; if $\alpha$ is simple then $\PP_\alpha$ will denote the corresponding standard parabolic subgroup and $\LL_\alpha$ a Levi subgroup of it. For any root $\alpha$ of $\AA$, $\check \alpha$ will denote the correponding co-root: $\Gm\to\AA$. We use additive, exponential notation for roots and co-roots -- for example, if $\check\alpha$ is a co-character into $\AA$ and $\chi$ is a character of $A$ then $e^{\check\alpha}(\chi)$ will denote $\chi(\check\alpha(\varpi))$. 

If $\YY$ is a $\BB$-variety with an open $\BB$-orbit, then $\mathring \YY$ denotes the open $\BB$-orbit. Given a $\BB$-variety $\YY$, we denote by $k(\YY)^{(\BB)}$ (resp.\ $k[\YY]^{(\BB)}$) the set of non-zero $\BB$-semiinvariants (eigenfunctions) on the rational (resp.\ regular) functions on $\YY$, and by $\varchi(\YY)$ the corresponding group of weights (eigencharacters). Finally, the space of an one-dimensional complex character $\chi$ of a group $H$ is denoted by $\CC_\chi$.

\begin{acknowledgements}{I would like to thank professors Joseph Bernstein, Daniel Bump, Dennis Gaitsgory, Herv\'e Jacquet, Friedrich Knop, Peter Sarnak and Akshay Venka-tesh for useful discussions, correspondence or references on this and related problems. Also, David Nadler for making available to me his preprints with Dennis Gaitsgory. Finally, the I.H.\'E.S. for its hospitality during the summer of 2007, when I completed the revision of this paper.}
\end{acknowledgements}


\section{Spherical varieties over algebraically closed fields} \label{secac}

\subsection{Basic notions} \label{ssspherical}

Let $\GG$ be an algebraic group over an arbitrary field $k$ in characteristic zero. By a $\GG$-variety (over $k$) we will mean a geometrically integral and separated $k$-scheme of finite type with an algebraic action of $\GG$ over $k$. A $\GG$-variety $\XX$ is called homogeneous if $\GG(\bar k)$ acts transitively on $\XX(\bar k)$ -- then $\XX$ is automatically non-singular. If $\XX$ has a point over $k$, its stabilizer $\GG_x$ is a subgroup over $k$ and $\XX\simeq\GG_x\backslash\GG$, the geometric quotient of $\GG$ by $\GG_x$. Conversely, for any closed subgroup $\HH$ the geometric quotient $\HH\backslash\GG$ is a homogeneous variety under the action of $\GG$.

Now, let $\GG$ be a reductive group over a field $k$. A normal $\GG$-variety $\XX$ over $k$ (not necessarily homogeneous) is called spherical if $\BB_{\bar k}$ (where $\BB_{\bar k}$ is a Borel subgroup of $\GG_{\bar k}$) has a Zariski open orbit on $\XX_{\bar k}$. This is equivalent (\cite{BrPr},\cite{Vi}) to the existence of finitely many $\BB_{\bar k}$-orbits. As a matter of convention, when we say ``a $\BB$-orbit on $\XX$'' (or ``$\GG$-orbit'') we will mean ``a $\BB_{\bar k}$-orbit on $\XX_{\bar k}$'' (respectively, ``$\GG_{\bar k}$''-orbit) -- then one naturally has to examine questions such as whether a ``$\BB$-orbit'' is defined over $k$, which will be the object of the next section.

For the whole paper, we will assume that $\XX$ is quasi-affine. This is not a really serious restriction: By \cite[Theorem 5.1]{Bo}, 
given a subgroup $\HH$ of $\GG$ there exists a finite-dimensional algebraic representation of $\GG$ over $k$ in which $\HH$ is the stabilizer of a line. If $\HH$ has trivial $k$-character group (i.e.\ group of homomorphisms $\HH\to\Gm$ over $k$), then this implies that $\HH\backslash\GG$ is embedded in the space of this representation and hence is quasi-affine. (Recall \cite[Proposition 1.8]{Bo} that an orbit of an algebraic group is always locally closed.) Hence, for an arbitrary $\HH$, we may replace $\HH$ by the kernel $\HH_0$ of all its $k$-characters and consider the quasi-affine variety $\HH_0\backslash \GG$, which is spherical for the $(\HH/\HH_0) \times \GG$ action.

From this point until the end of the present section we assume that $k$ is algebraically closed. Given a $\BB$-orbit $\YY$, the group of weights $\varchi(\YY)$ of $\BB$ acting on $k(\YY)$ is the character group of $\AA/\AA_Y$, where $\AA=\BB/\UU$ and $\AA_Y$ is the image modulo $\UU$ of the stabilizer of any point $y\in \YY$. The rank of $\varchi(\YY)$ is called the \emph{rank} of the orbit $\YY$. If $\YY$ is the open orbit, we will denote $\AA_Y$ by $\AA_X$\footnote{This will be a standard convention in our notation: if $\YY$ is a variety with an open $\BB$-orbit we will allow ourselves to use $\YY$ in the notation instead of $\mathring \YY$, whenever this causes no confusion.}; the corresponding rank is the \emph{rank} of the spherical variety. The rank of the open orbit is maximal among all $\BB$-orbits, as we explain below.

We recall the classification and properties of spherical subgroups $\HH$ for $\PPGL_2$.

\begin{theorem}[A classic] \label{classic} The spherical subgroups $\HH$ of $\GG=\PPGL_2$ over an algebraically closed field $k$ in characteristic zero are: 
\begin{description}
\item[Type G:] $\PPGL_2$; in this case there is a single $\BB$-orbit on $\HH\backslash \GG$.
\item[Type T:] a maximal torus $\TT$; there are three $\BB$-orbits: the open one, and two closed orbits of smaller rank.
\item[Type N:] $\mathcal{N}(\TT)$; there are two $\BB$-orbits: the open one, and a closed one of smaller rank.
\item[Type U:] $\SS\cdot \UU$, where $\UU$ is a maximal unipotent subgroup and $\SS\subset\mathcal N(\UU)$; there are two $\BB$-orbits, an open and a closed one, both of the same rank.
\end{description} 
\end{theorem}

\subsection{Knop's action} \label{ssaction} In \cite{KnOrbits}, F.~Knop defines an action of the Weyl group on the set of Borel orbits on a homogeneous spherical variety. This action is defined explicitly for simple reflections, and then it is shown that this description induces an action of the Weyl group (i.e. satisfies the braid relations). For the simple reflection $w_\alpha$ corresponding to a simple root $\alpha$ it is defined as follows:

Let $\PP_\alpha$ denote the parabolic (of semi-simple rank one) associated to $\alpha$ (for a fixed choice of maximal torus $\AA\hookrightarrow \BB$) and let $\YY$ be a $\BB$-orbit; the simple reflection $w_\alpha$ acts on the set of $\BB$-orbits contained in the $\PP_\alpha$-orbit $\YY\cdot \PP_\alpha$. Consider the quotient $\PP_\alpha \to \PPGL_2 = \Aut(\mathbb P^1)$ where $\mathbb P^1= \BB\backslash \PP_\alpha$. The image of the stabilizer $(\PP_\alpha)_y$ of a point $y\in\YY$ is a spherical subgroup of $\PPGL_2$, and according to the classification above we say that ``$(\YY,\alpha)$ is of type G, T, N or U''. (As a matter of language, we also say that ``$\alpha$ raises $\YY$ to $\ZZ$'' if $\ZZ\ne\YY$ is the open orbit in $\YY\PP_\alpha$.) We define the action according to the type of that spherical subgroup:

If it is of type G, $w_\alpha$ will stabilize the unique $\BB$-orbit in the given $\PP_\alpha$-orbit. In the case of type T, $w_\alpha$ stabilizes the open orbit and interchanges the other two. In the case of type N, $w_\alpha$ stabilizes both orbits. Finally, in the case of type U, $w_\alpha$ interchanges the two orbits. Since this defines a right action in our case that the group acts on the right, we modify it to a left action by defining $^w\YY := \YY^{w^{-1}}$, where $\YY^{w^{-1}}$ denotes the action of $w^{-1}$ on $\YY$ as defined by Knop; of course, in the case of simple reflections the description does not change. Notice that in every case the action of $w_\alpha$ preserves the rank of the orbit; more precisely, Knop proves:
\begin{lemma}\label{charlemma}
Let $\YY$ denote the open orbit of $\YY\PP_\alpha$, and let $\ZZ_*$ denote the closed orbits. There exist the following relations between their character groups: 
\begin{description}
\item[Type G:] $\varchi(\YY)\subset \left(\varchi(\AA)\right)^{w_\alpha}.$
\item[Type U:] ${^{w_\alpha}\!\varchi}(\ZZ) = \varchi(\YY)$.
\item[Type T:] ${^{w_\alpha}\!\varchi}(\ZZ_1) = \varchi(\ZZ_2)\subset \varchi(\YY)$.
\item[Type N:] ${^{w_\alpha}\!\varchi}(\ZZ) \subset \varchi(\YY)$.
\end{description}
\end{lemma}
(An exponent on the right denotes ``invariants''. An exponent on the left denotes the action of the Weyl group. Due to our modification of the definition, the lemma is true as stated, with the left action of $W$ on the characters.) In particular, $\varchi({^w\YY}) = {^w\!\varchi} (\YY)$ for every $w$ 
and the set $\mathfrak B_{00}$ of orbits of maximal rank is stable under the action of the Weyl group.

We denote the standard parabolic $\{g| \mathring \XX \cdot g = \mathring \XX\}$ (the elements of $\GG$ which preserve the open $\BB$-orbit) by $\PP(\XX)$. Equivalently, $\PP(\XX)$ is the parabolic corresponding to the simple roots $\alpha$ such that $\mathring \XX, \alpha$ is of type $G$. The Weyl group of its Levi will be denoted by $W_{P(X)}$. The little Weyl group $W_X\subset W$ of $\XX$ was mentioned, but not defined, in the introduction. The reader can take the following as the definition:

\begin{theorem}[Knop]
The stabilizer of $\mathring \XX$ under Knop's action is equal to $W_{(X)}:=W_X\ltimes W_{P(X)}$. The elements of $W_X$ are those of smallest length in their $W_{(X)}/W_{P(X)}$-coset.
\end{theorem}

\subsection{Parabolically induced spherical varieties}\label{ssparind}

There is an ``inductive'' process of constructing some spherical subgroups: Given a Levi subgroup $\LL$ of a parabolic $\PP\subset\GG$ and a spherical subgroup $\MM$ of $\LL$, we can form the subgroup $\HH=\MM\ltimes \UU_\PP$, which is a spherical subgroup of $\GG$. The structure of the $\BB$-orbits of $\XX:=\HH\backslash\GG$, relevant to the Borel orbits of $\MM\backslash \LL$, has been investigated by Brion \cite{BrOr}. The closure of each orbit $\YY$ of $\XX$ can be written uniquely as $\overline{\YY'w\BB}$ for $w\in [W_P\backslash W]$, where $[W_P\backslash W]$ denotes the set of representatives of minimal length for right cosets of $W_P$ (the Weyl group of $\LL$) and $\YY'$ a Borel orbit of $\XX':=\MM\backslash\LL$. We have $\varchi(\mathring\XX)=\varchi(\mathring\XX')$ and $W_X=W_{X'}$.

\subsection{Non-homogeneous spherical varieties}

Now we examine spherical varieties $\XX$ which are not necessarily homogeneous, i.e. may have more than one $\GG$-orbit. It is known then \cite{KnLV} that $\XX$ contains a finite number of $\GG$-orbits, and that each of them is also spherical. Let $\YY$ be a $\GG$-orbit. To $\YY$ one associates \cite[\S 2]{KnLV} the cone $\mathcal C_{\YY}(\XX)\subset \mathcal Q :=\Hom_{\mathbb Z}(\varchi(\mathring\XX),\mathbb Q)$ spanned by the valuations induced by $\BB$-stable prime divisors which contain $\YY$. This cone is non-trivial (more precisely \cite[Theorem 3.1]{KnLV}, there exists a bijection between isomorphism classes of ``simple embeddings'' of $\HH\backslash \GG$ and ``colored cones'') and we have: 

\begin{theorem}\label{extension} Let $\XX$ be a quasi-affine spherical variety, $\YY$ a $\GG$-orbit and $f\in k[\YY]^{(\BB)}$. There exists $f'\in k[\XX]^{(\BB)}$ with $f'|_{\YY}=f$. 
Hence, the group of weights of $\BB$ on $\mathring\YY$ is a subgroup of the weights of $\BB$ on $\mathring \XX$. 
More precisely, $\varchi(\YY)= \mathcal C_{\YY}(\XX)^\perp=\{\chi\in\varchi(\XX)| v(\chi)=0\textrm{ for every }v\in C_{\YY}(\XX)\}$. In particular, every non-open $\GG$-orbit on $\XX$ has strictly smaller rank than $\XX$ itself.
\end{theorem}

\begin{proof}
cf.\ \cite[Theorem 6.3]{KnLV}.
\end{proof}

\subsection{Non-degeneracy} \label{ssnondeg}

We recall the notion of a non-degenerate spherical variety \cite[\S 6]{KnOrbits},\cite{KnAs}:  The spherical variety $\XX$ is called non-degenerate if for every root $\alpha$ appearing in the unipotent radical of $\PP(\XX)$ there exists $\chi\in \varchi(\XX)$ such that $\chi^{\check\alpha}\ne 1$. This implies that $\PP(\XX)$ is the largest parabolic subgroup $\PP$ such that every character in $\varchi(\XX)$ extends to a character of $\PP(\XX)$. It is proven in \cite[Lemma 3.1]{KnAs} that every quasi-affine variety is non-degenerate. We will need a variant of this statement which includes the character groups of smaller $\BB$-orbits:

\begin{lemma}
Let $\XX$ be a quasi-affine spherical variety, and let $\YY$ be a $\BB$-orbit. Let $\alpha$ be a simple positive root that either does not raise $\YY$ (i.e. $\overline{\YY\PP_\alpha}=\overline{\YY}$) or raises $\YY$ of type U. Then either $(\YY,\alpha)$ is of type G (i.e. $\YY\PP_\alpha=\YY$) or there exists $\chi\in\varchi(\YY)$ with $\left<\chi,\check\alpha\right>\ne 0$.
\end{lemma}

\begin{proof}
Assume $\left<\chi,\check\alpha\right>=0$ for every $\chi\in\varchi(\YY)$. 
Recall that $\varchi(\YY)=\{\chi \,\, | \,\, \chi|_{\AA_Y}=1\}$; hence $\check\alpha(\Gm) \subset \AA_Y$. This cannot be the case if $\YY$ is the open orbit in $\YY\PP_\alpha$ and $(\YY,\alpha)$ is of type T or N. Let $(\YY,\alpha)$ be of type U -- without loss of generality, since $^{w_\alpha}\varchi(\YY)=\varchi(^{w_\alpha}\YY)$, $\YY$ is raised by $\alpha$. Given a point $y\in \YY$ with $\check\alpha(\Gm)\subset \BB_y$ (such a point must exist since all maximal tori of $\BB$ are conjugate inside of $\BB$), the Lie algebra of $\BB_y$ splits into a sum of eigenspaces of $\check\alpha(\Gm)$; if $\alpha$ raises $\YY$ of type U this implies that $\UU_\alpha \subset \BB_y$. Hence the stabilizer of $y$ in $[\LL_\alpha,\LL_\alpha]\simeq \SSL_2$ (or $\PPGL_2$) is a Borel subgroup $\BB_{\SSL_2}$ and we get an embedding of the complete variety $\BB_{\SSL_2}\backslash \SSL_2= \mathbb P^1$ into the quasi-affine variety $\XX$, a contradiction. Therefore $(\YY,\alpha)$ has to be of type G.
\end{proof}


\section{Rationality properties}\label{secrat}

\subsection{Homogeneous spaces}\label{ssgpsch}

The main questions that we examine in this section have to do with whether $\BB$- and $\GG$-orbits are defined over a non-algebraically closed field $k$ and to what extent Knop's action makes sense on the set of $k$-rational $\BB$-orbits. The results will be used in the next sections to examine the unramified spectrum over $p$-adic fields. We start by recalling certain classical results: 
We use the terminology of \cite{Bo}, according to which a solvable $k$-group is ``$k$-split'' (or simply ``split'') if it admits a normal series over $k$ whose successive quotients are $k$-isomorphic to $\Gm$ or $\Ga$ (in particular, connected).

\begin{theorem} \label{Borel}
Let $\GG$ be an algebraic group and $\HH$ a solvable algebraic subgroup. 
Assume that the maximal reductive quotient of $\HH$ is $k$-split.
\begin{enumerate}
\item \label{Borelpartone} If $\XX$ is a homogeneous $\HH$-variety then $\XX$ is affine and $\XX(k)\ne\emptyset$.
\item \label{Borelparttwo} $\GG(k)$ acts transitively on $(\HH\backslash\GG)(k)$.
\end{enumerate}
\end{theorem}

\begin{proof}
These are \cite[Theorem 15.11 and Corollary 15.7]{Bo}. Notice that in characteristic zero, every unipotent group is connected and $k$-split, therefore we only needed to assume that the quotient of $\HH$ by its unipotent radical was $k$-split in order to deduce that $\HH$ is $k$-split.
\end{proof}

\subsection{Rationality of the open Borel orbit} From now on, assume that $\GG$ is a split reductive group over a field $k$. This means that it has a Borel subgroup which is defined over $k$ and $k$-split. Let $\XX$ be a spherical $\GG$-variety (not necessarily homogeneous) over $k$. We assume that $\XX$ is quasi-affine (cf.\ \S \ref{ssspherical}).

As a generalization of Theorem \ref{Borel}.\ref{Borelpartone}, we prove:
\begin{proposition}\label{rationalorbit}
\begin{enumerate}
\item\label{rationalorbitpartone} Every line of $\BB$-eigenfunctions on $\bar k(\XX)$ is defined over $k$. 
\item\label{rationalorbitparttwo} The open $\BB$-orbit has a point (in particular, is defined) over $k$.
\end{enumerate}
\end{proposition}

The proposition is true in general for any quasi-affine variety $\XX$ over $k$ with a $k$-action of a split solvable group $\BB$ over $k$ such that $\BB$ has an open orbit on $\XX$.

\begin{proof}
The first claim follows from the fact that $\BB$ is split, hence all weights are defined over $k$, hence the (one-dimensional) eigenspaces for $\BB$ on $\bar{k}(\XX)$ are Galois invariant and, therefore, defined over $k$.

For the second, notice that there is a non-zero regular $\BB$-eigenfunction which vanishes on the complement of the open $\BB$-orbit. Indeed, the space of regular functions which vanish on the complement is non-zero (because $\XX$ is quasi-affine) and $\BB$-stable. As a representation of $\BB$ it decomposes into the direct sum of finite-dimensional ones. Let $V$ be such a finite-dimensional component. The space of $\UU$-invariants $V^\UU$ (where $\UU$ is the unipotent radical of $\BB$) is then non-zero, and since every finite-dimensional representation of $\AA=\BB/\UU$ is completely reducible, it follows that there exists a nonzero $\BB$-eigenfunction which vanishes on the complement of the open orbit.

Now, it follows from the first claim that this eigenfunction can be assumed to be in $k[\XX]$. Hence the open orbit is $k$-open (and, therefore, defined over $k$); by Theorem \ref{Borel}.\ref{Borelpartone}), the open $\BB$-orbit has a point over $k$.
\end{proof}

Because of this proposition, the open $\GG$-orbit on $\XX$ is isomorphic to $\HH\backslash\GG$ over $k$, where $\HH$ is a closed subgroup over $k$ (cf.\ \S \ref{ssgpsch}).

\subsection{Rationality of $\GG$-orbits}

With similar arguments as above we get:

\begin{proposition} \label{ratGorb}
 Every $\GG$-orbit on $\XX$ is defined over $k$. Hence, by Proposition \ref{rationalorbit}.\ref{rationalorbitparttwo} and the fact that all $\GG$-orbits are spherical, every $\GG$-orbit has a $k$-point on its open $\BB$-orbit.
\end{proposition}

\begin{proof}
Let $\ZZ$ be a $\GG$-orbit closure and let $\mathcal I\subset \bar k[\XX]$ be the ideal defining it. It is $\GG$-stable, hence a representation of $\GG$. It splits into a sum of irreducible finite-dimensional subspaces and every such is generated by a highest-weight vector, that is, a $\BB$-semiinvariant. By the rationality of those (Proposition \ref{rationalorbit}.\ref{rationalorbitpartone}), it follows that $\mathcal I$ is defined over $k$. Now by Proposition \ref{rationalorbit}, the open $\BB$-orbit, and hence the open $\GG$-orbit in $\ZZ$ are defined over $k$.
\end{proof}

\subsection{Splitting in {$\BB(k)$}-orbits}

Now we examine the splitting of the $k$-points of a $k$-rational $\BB$-orbit $\YY$ in $B=\BB(k)$-orbits. 

\begin{lemma}
 For every $k$-rational $\BB$-orbit $\YY$, the set of $B$-orbits on $Y$ is naturally a torsor for the (finite) abelian group $\Gamma_Y:=H^1(k,\AA_Y)$.
\end{lemma}

\begin{proof}
We know already that every $k$-rational $\BB$-orbit has a $k$-point $y$. It is known that for $k$-groups $\HH\subset\GG$ the $k$-orbits of $G$ on $(\HH\backslash\GG)(k)$ are parametrized (depending on the choice of an orbit) by the kernel of $H^1(k,\HH)\to H^1(k,\GG)$. Moreover, the first cohomology group of a unipotent group in characteristic zero is trivial, and so is, by Hilbert's Theorem 90, the first cohomology group of split tori. Hence $H^1(k,\BB)$ is trivial and $H^1(k,\BB_y)=H^1(k,\AA_Y)$. This proves the claim.
\end{proof}

Notice that by Lemma \ref{charlemma} and the fact that Knop's action is transitive on orbits of maximal rank, all $\AA_Y$, for $\YY$ of maximal rank, are $W$-conjugate to each other and, in particular, the order of $H^1(k,\AA_Y)$ is equal to $H^1(k,\AA_X)$ for all of them.

\begin{remark}
Again by Hilbert's Theorem 90, $H^1(k,\AA_Y)=H^1(k,\pi_0(\AA_Y))$, hence this group is non-trivial if and only if $\AA_Y$ is not connected. Notice also that $H^1(k,\AA_Y)$ has the following explicit description: it is equal to the quotient of $A'^Y$ by the image of $A$, where $\AA'^Y:=\AA/\AA_Y$. This can be seen from the long exact cohomology sequence for $1\to\AA_Y\to\AA\to\AA/\AA_Y\to 1$ and Hilbert's 90 again.
\end{remark}

\subsection{Spherical varieties for $\SSL_2$} \label{ssSL2}

Spherical varieties for $\SSL_2$ are of dimension at most 2. Therefore, homogeneous spherical varieties for $\SSL_2$ (over an arbitrary field $k$) belong to the homogeneous varieties classified by F.~Knop in \cite[Theorem 5.2]{KnR1}; it is easily seen that all the varieties in loc.\ cit.\ are spherical. We recall this classification according to the classification of the corresponding homogeneous spaces over the algebraic closure (Theorem \ref{classic}) and examine some basic rationality properties.

\subsubsection{Case G.}

The subgroup $\HH=\GG=\SSL_2$. There is a single $B$-orbit.

\subsubsection{Case T.}

The subgroup $\HH=\TT$, a (maximal) torus. 

By the equivalence of categories between diagonalizable $k$-groups and lattices with a Galois action, isomorphism classes of 1-dimensional tori over $k$ are classified by $\Hom(\Gal(\bar{k}/k),\mathbb Z/2)\simeq \widehat{(k^\times/(k^\times)^2)} $.

One way to describe the homogeneous space $\TT\backslash\GG$ is as $Q(1,\beta):= \{A\in\mathfrak{gl}_2 | \tr(A)=1, \det(A)=\beta\}$ where $\beta\in k$ and $4\beta-1\ne 0$, under the adjoint action of $\SSL_2$.

The space $X$ has, in general, several $G$-orbits. However, notice that we can naturally extend the action of $\SSL_2$ to an action of $\PPGL_2$, and $\PGL_2$ acts transitively on $X$.

We examine the splitting of the open Borel orbit in $B$-orbits: We have $\AA_X=\{\pm1\}$, and hence the orbits of $B$ on the open $\BB$-orbit are parametrized by $k^\times / {(k^\times)^2}$. (Don't confuse this with the parametrization of isomorphism classes of tori, mentioned above.) However, if we extend the action to $\PPGL_2$ then its Borel subgroup acts transitively on $\mathring X$ -- and this will be important later.

Notice that if the torus is non-split, it does not embed over $k$ into a Borel subgroup and therefore the smaller $\BB$-orbits do not have a point over $k$.

\subsubsection{Case N.}

The subgroup $\HH$ is the normalizer of a maximal torus. It turns out that for all tori we get the same homogeneous variety: Indeed, the space $\mathcal N(\TT)\backslash \SSL_2$ can be identified with the open subset of $\mathbb P(\rm{\mathfrak{pgl}_2})$ defined by $4\det(A)-(\tr(A))^2\ne 0$. It can then be seen that for every $\TT$, $\mathcal N(\TT)$ appears as a stabilizer of a $k$-point. Again, the action extends to $\PPGL_2$. Notice also that the $k$-points of $\mathcal N(\TT)$ coincide with the $k$-points of $\TT$ for $\TT$ non-split. This implies that the $\PGL_2$-orbit of a $k$-point with stabilizer $\mathcal N(\TT)$, for $\TT$ non-split, is isomorphic as a $\PGL_2$-space with $T\backslash\PGL_2$. The splitting of $\mathring X$ in $B$-orbits is parametrized by $k^\times/(k^\times)^4$, while if we consider the action of $\PPGL_2$ and let $\tilde\BB$ denote its Borel subgroup the orbits under $\tilde B$ are parametrized by $k^\times/(k^\times)^2$ (and the $B$-orbits are related to $\tilde B$-orbits through the natural map $k^\times/(k^\times)^4 \to k^\times/(k^\times)^2$).

\subsubsection{Case U.}

The subgroup $\HH$ is equal to $\SS\cdot \UU$, where $\UU$ is a maximal unipotent subgroup and $\SS\subset \mathcal N (\UU)$.

As a $G$-space, $X$ splits into a disjoint union of spaces isomorphic to $SU\backslash G$. The $k$-points of the open $\BB$-orbit may split into several $B$-orbits. However, because of the Bruhat decomposition over $k$, every one of them belongs to a different $G$-orbit. For the same reason, both $\BB$-orbits have $k$-points.

\subsection{Rationality of Knop's action}

\begin{proposition} \label{Knoprat}
\begin{enumerate} 
\item\label{Knopratpartone} If a $\BB$-orbit $\YY$ is defined (equivalently: has a point) over $k$, then so is ${^w\YY}$ for every $w\in W$. In particular, all the $\BB$-orbits of maximal rank are defined over $k$. 
\item\label{Knopratparttwo} If a $\BB$-orbit $\YY$ is raised by a simple root $\alpha$ to a $\BB$-orbit $\ZZ$, and $\YY$ is defined over $k$, then so is $\ZZ$. More precisely, if $y\in \YY(k)$ then $y\cdot P_\alpha$ contains a $k$-point of $\ZZ$.
\end{enumerate}
\end{proposition}

\begin{proof}
Consider the $\PP_\alpha$-orbit of $\YY$. Dividing by $\UU_{P_\alpha}$ we get a spherical variety for $\LL_\alpha$, the Levi of $\PP_\alpha$. Further dividing by the connected component of the center $\ZZ_\alpha$ of $\LL_\alpha$ we get a spherical variety $\XX_\alpha$ for $\SSL_2$ (or $\PPGL_2$). In both steps, the quotient maps are surjective on $k$-points since we are dividing by a unipotent group and a split torus, respectively. Therefore, a $\BB$-orbit on $\XX$ has a point over $k$ if its image in $\XX_\alpha$ is defined over $k$. By examining now the $\SSL_2$-spherical varieties which were classified above, it follows now that the rationality of $\YY$ implies the rationality of $^{w_\alpha}\YY$. The open orbit is, by Proposition \ref{rationalorbit}, defined over $k$. Moreover, all orbits of maximal rank belong to the open $\GG$-orbit and are in the $W$-orbit of the open orbit. This proves the rationality of all orbits of maximal rank. Finally, the fiber of the quotient map over the image of a point $z$ is acted upon transitively by $(\ZZ_\alpha^0/\ZZ_\alpha^0\cap \GG_z)(k)$ (these groups being canonically isomorphic for all $z\in(\YY\PP_\alpha)(k)$) therefore if $y'\cdot P_\alpha$ contains a point of $Z$, for $y'$ in the same fiber as $y$, then so does $y\cdot P_\alpha$.

\end{proof}

\subsection{The Zariski and Hausdorff topologies}

For any topological field $k$, the $k$-points of a variety $\XX$ over $k$ naturally inherit a topology from that of $k$. Indeed, since $X=\Hom_{k-\rm{alg}}(k[\XX],k)$ (we assume for simplicity of notation here that $\XX$ is affine, the general case can be recovered by covering $\XX$ by affine neighborhoods), every point can in particular be viewed as a map from $k[\XX]$ to $k$, and therefore the set of points inherits the compact-open topology from the space of such maps. (With $k$ having its given topology and $k[\XX]$ considered discrete.) If the topology on $k$ is locally compact, totally disconnected and Hausdorff, so will be the induced topology on $X$. We will conveniently refer to that topology as the ``Hausdorff'' topology.

For a spherical variety, we wish to examine the relation between closures of $B$-orbits in both the Zariski and the Hausdorff topology. By definition, essentially, the ``Hausdorff'' topology is finer than the Zariski topology, therefore a Zariski-open set is also Hausdorff-open. More precisely:

\begin{lemma}\label{open}
 If $\GG$ is a $k$-group acting on a $k$-variety $\XX$ with a (Zariski) open orbit, and if $x\in X$ belongs to the open orbit, then $x\cdot G$ is (Hausdorff) open in $X$.
\end{lemma}

\begin{proof}
 The differential $\mathfrak g\to T_x X$ is surjective, hence the claim.
\end{proof}

Does the Hausdorff closure of a $B$-orbit coincide with the Zariski closure? The example below shows that this is not the case, at least not in non-homogeneous varieties:

\begin{example} \label{ugly}Let $\XX$ be the subvariety of $\mathbb A^2\times (\mathbb A \smallsetminus \{0\})$ defined by the equation: $x^2-ay^2=0\,\,\,(a\ne 0)$. Consider the following action of $\GG=\BB=\Gm^2$: $(x,y,a)\cdot (r,k) = (rk^3x, rk^2y, k^2a).$
Then the $\BB(k)$ orbits are:
\begin{itemize}
\item $\{(x,y,\frac{x^2}{y^2})| x,y\ne 0\}$,
\item $\{(0,0,a)|a\in (k^\times)^2\}$,
\item $\{(0,0,a)|a\notin (k^\times)^2\}$.
\end{itemize}
The second and the third form together the $k$-points of the same $\BB$-orbit, but only the second is in the Hausdorff closure of the first.
\end{example}

Contrary to the previous example, for a homogeneous spherical variety we have:

\begin{lemma}\label{nbhd}
If $\XX$ is homogeneous then any neighborhood of a point $y\in Y$ (in the Hausdorff topology), where $\YY$ is a Borel orbit of dimension $j<\dim X$ contains $k$-points belonging to orbits of dimension $j+1$. In particular, every $G$-orbit contains points of $\mathring X$.
\end{lemma}

\begin{proof}
This is a direct consequence of Proposition \ref{Knoprat}.\ref{Knopratparttwo}: There is a simple root $\alpha$ raising $\YY$, and $y\cdot P_\alpha$ contains $k$-points of $^{w_\alpha}\YY$. By homogeneity, such points can be arbitrarily close to $y$.

The second claim now follows from the first and the fact that every $G$-orbit is open in the Hausdorff topology (Lemma \ref{open}).
\end{proof}

Our main object of study will be spaces of locally constant, compactly supported functions on spherical varieties. We need to decide whether, in the non-homogeneous case, we will allow the support of our functions to extend beyond the Hausdorff closure of $\mathring X$ (the $k$-points of the open $\BB$-orbit); or whether we will redefine $X$ as the Hausdorff closure of $\mathring X$. For the discussion of the next section and the results of section \ref{secKnop} it does not make a difference, since the smaller $\GG$-orbits have smaller rank and the results that we prove are ``generic'' and are not influenced by the smaller orbits. However, for the study of unramified vectors in section \ref{secEnd} we require that the support of all functions is contained in the Hausdorff closure of $\mathring X$. Of course, the Hausdorff closure of $\mathring X$ is $G$-stable. 

The above ``bad'' example will be understood better via the next lemma, which studies the relationship between $B$-orbits on the open $G$-orbit and open $B$-orbits on the smaller $G$-orbits.

\begin{lemma}\label{smorbitlemma}
For every $\GG$-orbit $\ZZ\subset\XX$ and every $k$-rational $\BB$-orbit closure $\VV$ containing $\ZZ$ there is a canonical $\AA$-equivariant map of geometric quotients: $\mathring\VV/\UU\to\mathring\ZZ/\UU$. Correspondingly, there is a canonical homomorphism of groups: $H^1(k,\AA_V)\to H^1(k,\AA_Z)$ such that the resulting map: \{$B$-orbits on $\mathring V$\} $\to$ \{$B$-orbits on  $\mathring Z$\} is equivariant. This map admits the following description: Each $B$-orbit on $\mathring V$ is mapped to the unique $B$-orbit on $\mathring Z$ which is contained in its closure -- in particular, the image of the map corresponds to the $B$-orbits of $\mathring Z$ which belong to the Hausdorff closure of $\mathring V$.
\end{lemma}

\begin{proof}
Recall (Theorem \ref{extension}) that every regular $\BB$-eigenfunction on $\ZZ$ extends to $\XX$ (in particular, to $\VV$). In other words, the canonical restriction $\varchi(\VV)\to\varchi(\ZZ)$ is surjective or, equivalently, $\AA_V\subset\AA_Z$. This induces:
$$\Gamma_V:=H^1(k,\AA_V)\to \Gamma_Z:=H^1(k,\AA_Z)$$ as claimed.

The extension property of eigenfunctions has the following consequence: The restriction: $k[\VV]\to k[\ZZ]$ splits canonically at the level of $\BB$-eigenfunctions:
$$ k[\ZZ]^{(\BB)}\hookrightarrow k[\VV]^{(\BB)}$$ or by passing to quotients:
$$ k[\mathring\ZZ]^{(\BB)}\hookrightarrow k[\mathring\VV]^{(\BB)}$$ which extends by linearity to $\UU$-invariants:
$$ k[\mathring\ZZ]^{\UU}\hookrightarrow k[\mathring\VV]^{\UU}.$$
Hence, we get a canonical morphism: $\mathring\VV/\UU\to\mathring\ZZ/\UU$, which of course is $\AA$-equivariant, and therefore induces a map between the sets of $B$-orbits which is compatible with the cohomology maps described above.

The implication: ``$zB\subset\overline{vB}$'' $\Rightarrow$ ``$vB\in \mathring V$ is mapped to $zB\in\mathring Y$''  follows immediately from the definition of the Hausdorff topology: Indeed, neighborhoods in this topology are determined by the values attained by regular functions, and if $vB$ is not mapped to $zB$ then this means that there exist $\BB$-semiinvariants strictly separating $vB$ from $zB$. To show the converse implication, assume that a neighborhood $N$ of $z$ does not meet a $vB$. Then the same is true for every $U$-translate of $N$, therefore we may assume that $N$ is $U$-invariant. But a fundamental system of $U$-invariant neighborhoods of $z$ is determined by the values of all $f\in k[\VV]^{(\BB)}$; therefore there exists a $B$-semiinvariant strictly separating $zB$ from $vB$ and therefore $vB$ is not mapped to $zB$.
\end{proof}

\subsection{Invariant differential forms and measures} \label{invariantmeasure}

Given a linear algebraic group $\GG$, its unipotent radical $\UU_G$ carries a (left and right) invariant top form $\omega$. It is unique up to scalar, and the adjoint (right) action: $\rm{Ad}_g:u\mapsto g^{-1}ug$ of $\GG$ transforms it by a character $\mathfrak d: \GG \to \Gm$ (the ``modular character''); in other words, $\rm{Ad}_g^*(\omega)= \mathfrak d(g) \omega$. This character is the sum of all roots of $\GG$ on the Lie algebra of $\UU_G$, and it is also equal to the ratio between a right- and a left-invariant top form on $\GG$ (which agree at the identity).

The group of isomorphism classes of $\GG$-line bundles on a homogeneous variety $\XX=\HH\backslash\GG$ (over the algebraic closure) is naturally: $\rm{Pic}_{\GG}(\XX):=\varchi(\HH)$ \cite{KKV}. Let $\mathcal L_\psi$ denote the corresponding line bundle for the character $\psi$ of $\HH$; its sections can be identified with global sections $f$ of the trivial bundle on $\GG$ such that $f(hg)=\psi(h)f(g)$ for all $h\in\HH, g\in\GG$. If $\psi$ is a $k$-character, then $\mathcal L_\psi$ is defined over $k$. There is a non-zero $\GG$-invariant global section of $\mathcal L_\psi^*\otimes \Omega$ (the sheaf of top-degree differential forms valued in the dual of $\mathcal L_\psi$) if and only if $\psi=\frac{\mathfrak d_\HH}{\mathfrak d_\GG|_\HH}$; in particular, there is a $\GG$-invariant top form on $\XX$ if and only if the modular characters of $\GG$ and of $\HH$ agree on $\HH$. 

Given a smooth variety $\XX$ over a local field $k$, any $k$-rational top differential form $\omega$ on $\XX$ gives rise to a positive Borel measure on the topological space of its $k$-rational points \cite{We}. This measure will be denoted by $|\omega|$. The complex character $\delta_G:=|\mathfrak d_\GG|$ of $G=\GG(k)$ is equal to the ratio between right and left Haar measure on $G$ and is also called ``the modular character''. 

We will show that without loss of generality we may assume the $k$-points of the open orbit $\mathring\XX$ in our spherical variety $\XX$ to possess a $B$-invariant measure. For this, we may take $\XX$ to be homogeneous.

As discussed in \S \ref{ssspherical}, we can assume that $\GG=\GG_1\times\TT$ over $k$, where $\TT$ is a torus, $\GG_1$ acts transitively on $\XX$: $\XX=\HH_0\backslash \GG_1$, where $\HH_0$ has no $k$-characters. Then $\XX$ possesses a $\GG_1$-invariant $k$-rational top form $\omega$. The idea is to replace $\TT$, if necessary, by a subtorus.

Let $\BB_1=\GG_1\cap \BB$, a Borel subgroup of $\GG_1$; then $\BB=\BB_1\times\TT$. For the open orbit we have: $\mathring\XX=(\HH\cap (\BB_1\times\TT))\backslash (\BB_1\times\TT)$ (assuming that $\HH\BB$ is open). Then the quotient $\mathring\XX/\BB_1\simeq (\HH\BB_1\cap\TT)\backslash \TT$. Let $\HH_\TT=\HH\BB_1\cap\TT$. Then $\TT$ is a finite quotient over $k$ of $\TT'\times\HH_T$, where $\TT'$ is some subtorus and $\XX$ is still a spherical $\GG_1\times \TT'$-variety. Let $\chi$ be the $k$-character under which $\omega$ transforms under the action of $\TT$. Since $\HH_\TT\cap\TT'$ is finite, it follows that $|\omega|$ (which is a positive measure on $X$) is invariant under $(\HH_\TT\cap\TT')(k)$. Hence $|\omega|$ varies by a positive (unramified) character of $(\TT_H\cap\TT')(k)\backslash \TT'(k)$, and by twisting it by the inverse of that character (which is constant on the orbits of $\BB_1$) we obtain a $B_1\times T'$-invariant measure on $\mathring X$.


\section{Mackey theory and intertwining operators}\label{secMackey}

In this section we summarize the method of intertwining operators (cf.\ \cite{BZ,Cas}). It is usually referred to as ``Mackey theory'' by analogy to Mackey's theorem for representations of finite groups. The finiteness of $\BB$-orbits is very important here. We use a method of Igusa to establish the rationality and other properties of intertwining operators, we examine their poles and we discuss the precise relationship between intertwining operators constructed analytically and Jacquet modules.

\subsection{Unramified principal series} \label{ssups}

From now on, $k$ will always denote a locally compact non-archimedean local field in characteristic zero.
We work in the abelian category $\mathcal S$ of smooth representations of $G$, which means that every vector has open stabilizer. 

An \emph{unramified character} of a reductive algebraic group $A$ over a $p$-adic field $k$ is a complex character of $A$ which is of the form $|f_1|^{s_1}\cdots |f_r|^{s_r}$, where $f_1,\dots,f_r$ are $k$-rational algebraic characters of $\AA$ (i.e.\ homomorphisms into $\Gm$), defined over $k$, the sign $|\cdot|$ denotes the $p$-adic absolute value and $s_1,\dots, s_r \in \CC$.

The group of unramified characters of $A$ a natural structure of a complex algebraic torus: If $f_1,\dots,f_m$ form a basis for the group of algebraic characters modulo torsion and if $\chi=|f_1|^{s_1}\cdots |f_m|^{s_m}$ then the association $\chi\mapsto (q^{-s_1},\dots,q^{-s_m})\in (\CC^\times)^m$, where $q$ is the order of the residue field of $k$, defines the structure of a complex torus on the group of unramified characters.

If $\XX=\HH\backslash\GG$ is a homogeneous $\GG$-variety over $k$, an unramified character $\psi$ of $H$ gives rise to a complex $G$-line bundle over $X$, to be denoted by $\mathcal L_\psi$. If $X=H\backslash G$ then sections of this line bundle can be described as complex functions $f$ on $G$ such that $f(hg)=\psi(h)f(g)$ for every $h\in H, g\in G$. In general, choose a normal subgroup $\HH_1\subset \HH$ such that $H_1$ is in the kernel of all unramified characters of $H$ and $\HH/\HH_1$ is a $k$-split torus (this is always possible: choose a quotient $\HH/\HH_1$ of $\HH/\HH_0$ -- where $\HH_0$ is as in \S\ref{ssspherical} -- such that $\varchi(\HH/\HH_1)_k$ is isomorphic to the quotient of $\varchi(\HH/\HH_0)_k$ by its torsion) then sections of $\mathcal L_\psi$ can be described as functions $f$ on $(\HH_1\backslash\GG)(k)$ such that $f(hx)=\psi(h)f(x)$ for $h\in (\HH/\HH_1)(k)$. (Recall that the quotient map $\HH_1\backslash\GG\to\HH\backslash\GG$ is surjective on $k$-points.) There is an $\mathcal L_\psi^*$-valued $G$-invariant measure on $X$ if and only if $\psi=\frac{\delta_H}{\delta_G}$.

Let $B$ be the Borel subgroup of $G$, with a maximal torus $A\subset B$; we will denote the complex torus of its unramified characters (considered simultaneously as characters of $B$ via $A=B/U$) by $A^*$. This is the maximal torus in the ``Langlands dual'' group of $\GG$. Co-roots of $\AA$ are naturally roots of $A^*$ and hence the expression $e^{\check\alpha}(\chi)$ ($\chi\in A^*$) makes sense and is equal to $\chi(e^{\check\alpha}(\varpi))$. (This is compatible with the standard conventions of \cite{BoAut}.) Let $\mathfrak d$ (resp.\ $\delta$) be the algebraic (resp.\ complex) modular character of the Borel; hence $\mathfrak d = e^{2\rho}$ (where $\rho$ is the half-sum of positive roots) and $\delta=|\mathfrak d|$. Given an unramified character $\chi$ of $B$, we define the unramified principal series $I(\chi):=\Ind_B^G(\chi\delta^\frac{1}{2})$. (Since we are working in the smooth category, $\Ind_B^G(\chi\delta^\frac{1}{2})$ is the space of \emph{smooth} sections of $\mathcal L_{\chi\delta^\frac{1}{2}}$ over $B\backslash G$.) We recall its properties: For a hyperspecial maximal compact subgroup $K$ of $G$, it contains a unique (up to scalar multiple) vector invariant under $K$ (called ``unramified''). For generic $\chi\in A^*$, it is irreducible. Generic, when talking about points on complex varieties, will mean ``everywhere, except possibly on a finite number of divisors''. For generic $\chi$, again, and $w\in W$ we have an isomorphism $T_w: I(\chi)\simeq I({^w\chi})$. We will recall the construction of the intertwining operator $T_w$ later. Also, the spaces $I(\chi)$ can be identified to each other as vector spaces by considering the restriction of $f\in I(\chi)$ to $K$. If we call this common underlying vector space $V$, and we have a family of maps $m_\chi$ from a set $S$ to $I(\chi)$ for $\chi$ varying on a subvariety $D$ of $A^*$, then we say that the family is \emph{regular} if for every $s\in S$ we have $m_\chi(s)\in V\otimes \CC[D]$. Similarly we define the notion of ``rational'' family of maps. We write $m_{\chi,1}\sim m_{\chi,2}$ to denote that $m_{\chi,1}=c(\chi) m_{\chi,2}$ for some non-zero rational function $c$ of $\chi$.

We will need to recall more information on the divisors on which the above statements (irreducibility of $I(\chi)$ and isomorphism with $I({^w\chi})$) may fail to be true: First, there are the ``irregular'' characters, i.e.\ those given by an equation $\chi={^w\chi}$, $w\in W$. Those are precisely the characters belonging to one of the divisors $R_{\check\alpha}:=\{\chi|\chi^{\check\alpha}=1\}$, where $\check\alpha$ is a co-root. More precisely, the representation $I(\chi)$ \emph{may} be reducible for $\chi$ irregular, and the intertwining operator $T_w$ has a pole on the divisor $\bigcup_{\check\alpha>0, w\check\alpha<0} R_{\check\alpha}$. Then, there are the divisors $Q_{\check\alpha}$ ($\check\alpha$ a co-root) described by the equation $\chi^{\check\alpha}=q$. It is known that, for such $\chi$, $I(\chi)$ is reducible and $T_w$ ceases to be an isomorphism on the divisor:
$$ \bigcup_{\check\alpha>0,w \check\alpha<0} (Q_{\check\alpha} \cup Q_{-\check\alpha}). $$

Returning to our spherical variety, the complex torus of unramified characters of $B$ supported by an orbit $Y$ (i.e.\ generated by complex powers of the modulus of $k$-rational $\BB$-semiinvariants on $\YY$) will be denoted by $A_Y^*$. If $\YY=$ the open $\BB$-orbit on $\XX$ we will be denoting $A_Y^*$ by $A_X^*$. (Its Lie algebra is $\mathfrak a_X^*=\varchi(\XX)\otimes_{\mathbb Z}\CC$.)

\subsection{The Bernstein decomposition and centre}\label{ssBernstein} By the theory of the Bernstein centre \cite{BeCtr}, the category $\mathcal S$ is the direct sum of categories $\mathcal S_{P,\sigma}$, indexed by \emph{equivalence classes} of pairs of data (``parabolic subgroup'', ``orbit of irreducible supercuspidal representations of its Levi''). (Here ``orbit'' implies the action of the torus of unramified characters of the Levi subgroup and two such sets of data are equivalent if and only if they are conjugate by some $g\in G$.) The ``simplest'' of these categories is indexed by the data (``Borel subgroup'', ``unramified characters''). It will, by abuse of language, be called the ``unramified Bernstein component'', although not all representations belonging to it are unramified (i.e.\ possess a vector invariant under a maximal compact subgroup). Given a smooth representation $\pi$, its ``unramified'' direct summand $\pi_\ur$ admits the following equivalent characterizations:
\begin{enumerate}
\item Every irreducible subquotient of $\pi_\ur$ and no irreducible subquotient of its complement is isomorphic to a subquotient of some unramified principal series.
\item $\pi_\ur$ is the space generated by the vectors of $\pi$ which are invariant under the Iwahori subgroup.
\end{enumerate}

Moreover, the centre $\mathfrak z(\mathcal S)$ of $\mathcal S$ is described in \cite{BeCtr}: This is, by definition, the endomorphism ring of the identity functor; in other words, every element of this ring is a collection of endomorphisms, one for each object in the category, such that when applied simultaneously they commute with all morphisms in the category. The centre can also be identified with the convolution ring of all conjugation-invariant distributions on $G$ whose support becomes compact when they are convolved with the characteristic measure of any open-compact subgroup.

By the above decomposition, one evidently has $\mathfrak z(\mathcal S) = \prod_{P,\sigma} \mathfrak z(\mathcal S_{P,\sigma})$.
Each of the factors in this product is naturally isomorphic to the space of regular functions on a complex variety (and the disjoint union of these varieties is called the ``Bernstein variety''). The centre of $\mathcal S_\ur$ is naturally isomorphic to $\CC[A^*]^W$ by mapping each element to the scalar by which it acts on $I(\chi)$, for all $\chi\in A^*$. Convolving the corresponding distributions with the characteristic measure of a hyperspecial maximal compact subgroup $K$, we get an isomorphism of rings between $\mathfrak z(\mathcal S_\ur)$ and the ``spherical Hecke algebra'' $\mathcal H (G,K)$ of $K$-biinvariant measures on $G$. The fact that $\mathcal H (G,K)\simeq \CC[A^*]^W$ is the \emph{Satake isomorphism}.

\subsection{Filtrations} Let $X$ be a locally compact, totally disconnected space with a continuous (right) action of $G$. Then the space $C_c^\infty(X)$ of locally constant, compactly supported complex functions on $X$ furnishes a smooth representation of $G$, the ``right regular representation'', to be denoted by $g\mapsto R(g)$. The discussion below applies more generally to the space $C_c^\infty(X,\mathcal L_\psi)$ of smooth, compactly supported sections of $\mathcal L_\psi$, where $\mathcal L_\psi$ is as in \S\ref{ssups}, but for simplicity we will work with the trivial bundle here and make a few comments on the general case in \S \ref{ssbundle}.

By Frobenius reciprocity,\footnote{The isomorphism asserted by Frobenius reciprocity is given as follows: Given a morphism into $I(\chi)$, compose with ``evaluation at 1'' to get a functional into $\CC_{\chi\delta^{\frac{1}{2}}}$. For the whole paper, we will avoid distinguishing between the morphism and the functional whenever possible, and we will be using the same letter to denote both.} $$\Hom_G \left(C_c^\infty(X), I(\chi)\right) = \Hom_B \left(C_c^\infty(X), \CC_{\chi\delta^{\frac{1}{2}}}\right).$$ If $Y\subset X$ is open and $B$-stable, and $Z= X\smallsetminus Y$, then we have an exact sequence:
$$ 0\to C_c^\infty(Y) \to C_c^\infty(X) \to C_c^\infty(Z) \to 0$$
which gives rise to an exact sequence of distributions (by definition, the linear dual of $C_c^\infty$):
$$ 0\to \mathcal D(Z) \to \mathcal D(X) \to \mathcal D(Y) \to 0.$$
By applying the functor of $(B,\chi^{-1}\delta^{-\frac{1}{2}})$-equivariance we get a sequence on the spaces $\mathcal D(\bullet)^{(B,\chi^{-1}\delta^{-\frac{1}{2}})} \linebreak = \Hom_B \left(C_c^\infty(\bullet), \CC_{\chi\delta^{\frac{1}{2}}}\right)$ (recall that, by definition, the action of $g$ on a distribution $D$ is given by $\pi^*(g)D (f) = D (\pi(g^{-1})f)$), but we might lose right exactness:
\begin{equation}\label{equivdistr}
 0\to\mathcal D(Z)^{(B,\chi^{-1}\delta^{-\frac{1}{2}})} \to \mathcal D(X)^{(B,\chi^{-1}\delta^{-\frac{1}{2}})} \to \mathcal D(Y)^{(B,\chi^{-1}\delta^{-\frac{1}{2}})}.
\end{equation}

We apply the above in the setting of $X$= the $k$-points of our spherical variety $\XX$. As we remarked, the Zariski topology is coarser than the induced Hausdorff topology, hence the set of $k$-points of orbits of dimension $>d$ is open in the set of $k$-points of orbits of dimension $\ge d$. More precisely, we have the following filtration:
$$ 0 \hookrightarrow C_c^\infty(\mathring X) \hookrightarrow C_c^\infty(\cup_{\dim(\YY)\ge\dim(\XX)-1} Y) \hookrightarrow \dots \hookrightarrow C_c^\infty(X)$$
with successive quotients isomorphic to $C_c^\infty(\cup_{\dim(\YY)=d} Y)$ for the appropriate $d$.

It follows that the dimension of the space of $(B,\chi^{-1}\delta^{-\frac{1}{2}})$-equivariant distributions on $X$ is less than or equal to the sum of the dimensions of $(B,\chi^{-1}\delta^{-\frac{1}{2}})$-equivariant distributions on $Y$, for all $k$-rational orbits $\YY$. Part of what we prove below is that, for generic $\chi$, we actually have equality. In any case, the problem now has been divided in two parts: examine the question of $(B,\chi^{-1}\delta^{-\frac{1}{2}})$-equivariant distributions on any single $\BB$-orbit, and then determine whether these distributions extend to the whole space.

\subsection{Distributions on a single orbit}
\subsubsection{The case where $B$ acts transitively.} Now we examine the question of \linebreak
$(B,\chi^{-1}\delta^{-\frac{1}{2}})$-equivari-ant distributions on $C_c^\infty(Y)$, where $Y=y \cdot B$ is some orbit of $B$ (not of $\BB$). We will use the natural projection: $C_c^\infty(B)\twoheadrightarrow C_c^\infty(y B)= C_c^\infty(B_y\backslash B)$, given by integration on the left over $B_y$ with respect to right Haar measure:
\[ p(f)(x)=\int_{B_y} f(bx) d_r b \,\,\, ,f\in C_c^\infty(B) \] in order to pull-back such a distribution to $B$.
If we pull back a $(B,\chi^{-1}\delta^{-\frac{1}{2}})$-equivariant distribution to $B$, we get a $(B,\chi^{-1}\delta^{-\frac{1}{2}})$-equivariant distribution on $B$, which has to be $\chi^{-1}\delta^{-\frac{1}{2}}$ times right Haar measure, hence equal to $\chi^{-1}\delta^{\frac{1}{2}}\cdot d_l b$, where $d_l b$ denotes \emph{left} Haar measure. In other words the distribution is given by 
\begin{equation} \label{SY1}
S_\chi^Y:\phi\mapsto \int_B f(b) \chi^{-1}\delta^\frac{1}{2}(b) d_l b
\end{equation}
where $f\in C_c^\infty (B)$ such that $p(f)=\phi$. This distribution is well defined (i.e.\ will factor through the surjection $p: C_c^\infty(B) \twoheadrightarrow C_c^\infty (B_y \backslash B)$) if and only if: 
\begin{equation} \label{condition}
\left.\chi^{-1}\delta^\frac{1}{2} \right|_{B_y}=\delta_{B_y}
\end{equation}
where $\delta_{B_y}$ is the modular character of $B_y$.

\begin{definition}
Given a $B$-orbit $Y$, the unramified characters satisfying (\ref{condition}) will be called $Y$-admissible. They form a complex subvariety of $A^*$ which will be denoted by $\Adm_Y$.
\end{definition}

As a matter of notation, for $\YY$ a $\BB$-orbit closure we will denote $\Adm_Y=\Adm_{\mathring Y}$. Notice that two $Y$-admissible characters $\chi$ differ by an element of $A_Y^*$. Hence, \emph{$\Adm_Y$ is a translate of $A_Y^*$; in particular, $\dim\Adm_Y=\rk(Y)$}. Clearly, $\Adm_Y$ only depends on the $\BB$-orbit containing $Y$. Moreover, it is immediate that \emph{the family of functionals $S_\chi^Y$ is regular in $\chi\in\Adm_Y$ (cf.\ \S \ref{ssups}).}

Summarizing the above discussion:

\begin{lemma}
For each $B$-orbit $Y$ we have a (unique up to scalar) $B$-morphism: $C_c^\infty(Y)\to \CC_{\chi\delta^{\frac{1}{2}}}$ if and only $\chi\in\Adm_Y$. These morphisms are given as specializations of a regular family $S_\chi^Y$, which is uniquely defined up to a non-vanishing regular function on $\Adm_Y$.
\end{lemma}

(Recall that a non-vanishing regular function on a group variety is always the multiple of a character.)

\begin{remark}
All the rational families of distributions/intertwining operators to be defined in this paper are uniquely defined up to a non-vanishing regular function on the corresponding parametrizing variety. This dependence is typically originating from a choice of base point which is used in order to write down certain integral expressions, and choices of measures etc. The relations to be established between such distributions are always of the form $S_1\sim S_2$ (cf.\ \S\ref{ssups}), and whenever we write $S_1=S_2$ it should be interpreted as ``equality up to an invertible regular function of $\chi$'', which is slightly stronger than $S_1\sim S_2$. In fact, for the results of this paper the reader might as well ignore the normalization up to a regular function, and consider the morphisms as if they were uniquely defined up to a rational function.
\end{remark}

In case that the orbit $Y$ admits a $B$-invariant measure, or equivalently that $\delta|_{B_y}=\delta_{B_y}$, the condition of admissibility takes the nicer form: $$\chi\delta^\frac{1}{2}\in A_Y^*$$ and the distribution can be expressed as an integral on $Y$: 
\begin{equation} \label{SY2}
S_\chi^Y:\phi\mapsto \int_Y \phi(y) |f_1(y)|^{s_1}\cdots|f_m(y)|^{s_m} dy
\end{equation}
where $f_1,\dots, f_m$ are a basis for the $k$-rational semiinvariants of $\BB$ on $\YY$ (modulo torsion) and $s_1,\dots,s_m\in \CC$ such that $|f_1|^{s_1}\cdots |f_m|^{s_m}$ is of weight $\chi^{-1}\delta^{-\frac{1}{2}}$. More generally, even if there is no $B$-invariant measure there will always be a measure $dy$ that varies by some unramified character $\psi$ of $B$ (since every unramified character of $B_y$ can be extended to a character of $B$) and then the same expression will give $S_\chi^Y$ except that the weight of $|f_1|^{s_1}\cdots |f_m|^{s_m}$ should be $\psi\cdot\chi^{-1}\delta^{-\frac{1}{2}}$.

\subsubsection{Non-transitive action of $B$ and weighted distributions.} \label{ssweighted} Let $\YY$ be a $k$-rational $\BB$-orbit and denote by $\underline S_\chi^Y$ the vector space of $B$-morphisms: $C_c^\infty(Y)\to \CC_{\chi\delta^\frac{1}{2}}$. It admits a basis consisting of $S_\chi^{Y_i}$, where $Y_i$, $i=1,\dots,n$ ranges over the distinct $B$-orbits; however, this turns out not to be the correct basis to use. In this section we define certain bases which will become useful later.

The first basis which we will defined is quite natural: Fix a point $y\in Y$ thus getting an identification: $\YY/\UU\simeq \AA'^Y:=\AA/\AA_Y$. Notice that we can choose $\AA\subset\BB$ such that $\AA_\YY\subset\BB_y$, hence making the torus $\AA'^Y$ act on $\YY$ ``on the left'' -- the resulting action depends on the choices made, but the choices will not have any effect to our definitions except (as usual) up to non-vanishing regular functions of $\chi$.  We may form a basis for $\underline S_\chi^Y$ indexed by the set of complex characters $\tilde\chi$ of $A'^Y$ which coincide with $\chi$ on the image of $A$.\footnote{There is a slight abuse of language here since $\chi$ is not a character of $A/A_Y$ but, rather, satisfies equation (\ref{condition}). For convenience, we will say that ``$\tilde\chi$ is a character of $A'^Y$ which extends $\chi$'' to mean that $\tilde\chi$ must also satisfy (\ref{condition}).} The corresponding basis element $S_{\tilde\chi}$ is defined in a similar way as in the previous section. To be precise, $S_{\tilde\chi}^Y$ can be described as the composition of two morphisms. The first is integration over the horocycles on $Y$ (with some abuse of notation \footnote{The abuse has to do with the fact that $A'^Y$ does not act on the right. The proper interpretation of the integral is as follows: First, recall the definition of $\mathcal L_{\delta\delta_{B_y}^{-1}}$ from \S\ref{ssups}. If $\AA_1$= the kernel of the algebraic character $\mathfrak d \mathfrak d_{B_y}^{-1}$ in $\AA_Y$ then $f$ is a function in $\Ind_{(\AA_Y/\AA_1)(k)}^{(\AA/\AA_1)(k)}(\delta\delta_{B_y}^{-1})$. Secondly, for every $a \in (\AA_Y/\AA_1)(k)$ pick a preimage $\bar a\in\AA(\bar k)$; then conjugation by $\bar a$ carries a $\UU$-invariant top form on $\UU_y\backslash\UU$ to a $\UU$-invariant top form on $\UU_{y\cdot a}\backslash\UU$, necessarily $k$-rational. These forms give rise to measures on the sets of $k$-points, and with respect to those measures: $f(a)=\int_{U_{y\cdot a}\backslash U} \phi(y\cdot a \cdot u) du$.} and depending, again, up to non-zero rational functions of $\tilde\chi$, on the choices of invariant measure and of base point):

\begin{equation}\label{Uint}
 C_c^\infty(Y) \to C_c^\infty(A'^Y, \mathcal L_{\delta\delta_{B_y}^{-1}}): f(\bullet)= \int_{U_y\backslash U} \phi(y\cdot u\bullet) du
\end{equation}
followed by integration over the torus:
\begin{equation}\label{inttorus} \int_{A'^Y} f(a) \tilde\chi^{-1}\delta^{-\frac{1}{2}}(a) da. \end{equation}
(Here notice that by the admissibility condition (\ref{condition}) the product $f(a) \tilde\chi^{-1}\delta^{-\frac{1}{2}}(a)$ lies in the trivial line bundle over $A'^Y$.)

While this basis is very natural, unfortunately in certain cases one needs to use yet a different one: While the $S_\chi^{Y_i}$'s are distributions supported on a single $B$-orbit, and the $S_{\tilde\chi}^Y$'s are weighted averages of the former as they range over \emph{all} $B$-orbits on $Y$, the new basis will be somewhere between the two: It will consist of weighted averages over \emph{some} of the $S_\chi^{Y_i}$'s. Since we will only use such a basis in a very specific case, we give the definition only for that case. At this point, the definition will appear very unmotivated, and the reader should skip it at first reading. 

The definition will depend not only on $Y$ but also on some simple root $\alpha$ such that $(Y,\alpha)$ is of type N with $Y$ being the open orbit in $YP_\alpha$.
We want to define first a subgroup $A'^Y_\alpha$ of $A'^Y$. We consider the quotient map: $\YY\PP_\alpha \to \mathcal N(\TT)\backslash \PPGL_2$ and let $\BB_2$ denote the corresponding Borel subgroup of $\PPGL_2$. Then we have a map from $B$-orbits on $Y$ to $B_2$-orbits on the image of $Y$, corresponding to the map 
\begin{equation}\label{chmap}
H^1(k,\AA_Y)\to H^1(k,\AA_{Y,2})
\end{equation}
where $\AA_{Y,2}$ is the image of $\AA_Y$ modulo the center of $\LL_\alpha$. Recall that $H^1(k,\AA_Y)$ is a quotient of $A'^Y$; now we define $A'^Y_\alpha$ to be the preimage in $A'^Y$ of the kernel of (\ref{chmap}).

Now for every character $\tilde\chi$ of $A'^Y_\alpha$ extending $\chi$, and for every orbit $\zeta$ of $A'^Y_\alpha$ on $Y/U$ (the latter being naturally a torsor of $A'^Y$) we can define a morphism $S_{\tilde\chi,\zeta}^Y$ as above, except that we restrict the last integration (\ref{inttorus}) to the chosen orbit of $S_{\tilde\chi,\zeta}^Y$. In other words, we take $y\in Y$ with $(y\mod U)\in \zeta$ in order to define the identification $\YY/\UU\simeq \AA'^Y$ and repeat the first step (\ref{Uint}) while replacing (\ref{inttorus}) by:
\begin{equation}\label{intsmalltorus} \int_{A'^Y_\alpha} f(a) \tilde\chi^{-1}\delta^{-\frac{1}{2}}(a) da. \end{equation}

Again, our regular family of morphisms is only well-defined up to a non-vanishing regular function. The new basis for $\underline S_\chi^Y$ is indexed by pairs $(\tilde\chi,\zeta)$ where $\tilde\chi$ is an extension of $\chi$ to $A'^Y_\alpha$ and $\zeta$ denotes an orbit of $A'^Y_\alpha$ on $Y/U$. Notice that in the notation we suppress the dependence of this basis definition on $\alpha$; however, this dependence is certainly an unpleasant feature which complicates the final results.

\subsubsection{Comparison between admissible characters for different orbits.} \label{sssadmissible} Based on \linebreak Lemma \ref{charlemma}, we can describe the relations between the varieties of admissible characters on the Borel orbits of a $\PP_\alpha$-orbit.

\begin{description}
\item[Case G]: We have $\chi^{\check\alpha}=q$ for every $\chi\in\Adm_Y$. (Recall that $\alpha$ denotes the positive root in the Levi of $\PP_\alpha$.)
\item[Case U]: Notice that the stabilizer $\BB_z$ of the closed orbit has unipotent radical of dimension one larger than the stabilizer $\BB_y$ of the open orbit. In fact, $\delta^{\frac{1}{2}}|_{B_z}= e^\alpha\cdot {^{w_\alpha}\delta^\frac{1}{2}_{B_y}}$ and $\delta_{B_z}={^{w_\alpha}\delta_{B_y}}$ (as characters on $A_Z={^{w_\alpha}}\!A_Y$), and this implies that $\Adm_{^{w_\alpha}\!Y}={^{w_\alpha}\!\Adm_Y}$.
\item[Case T:] If $\YY$ denotes the open orbit and $\ZZ_*$ the closed ones (in the case of a split torus, for we have seen that if $\TT$ is non split then the closed orbits are not defined over $k$) then $\Adm_Y\supset\Adm_{Z_*}$. Notice, however, that for the small orbits it does not hold that $\Adm_{Z_1}= {^{w_\alpha}\!\Adm_{Z_2}}$ -- the correct relation is $\delta^{-\frac{1}{2}}\Adm_{Z_1}={^{w_\alpha}\!\!\left( \delta^{-\frac{1}{2}}\Adm_{Z_2}\right)}$.
\item[Case N:] If $\YY$ denotes the open orbit and $\ZZ$ the closed one then $\Adm_Y\supset\Adm_Z$ and $\delta^{-\frac{1}{2}}\Adm_Z$ is $w_\alpha$-invariant.
\end{description}

\begin{remark}
 Notice that in case $T$ we do not have $^{w_\alpha}(\Adm_{Z_1})\subset \Adm_{Z_2}$. (Similarly, in case $N$ we do not have this for $Z_1=Z_2=Z$). Moreover, it may be contained in $\mathcal Q_{-\alpha}$ but not in $\mathcal Q_\alpha$ or $\mathcal R_\alpha$. (For the definition of $\mathcal Q_{\alpha}$ and $\mathcal R_\alpha$ cf.\ \S\ref{ssups}.) Indeed, the condition $\chi^{-1}\delta^\frac{1}{2}|_{B_z}=\delta_{B_z}$ and the fact that $B_{z_1}={^w B_{z_2}}$ imply that:
\begin{enumerate}
\item Either $A^*_{Z_1}=A^*_{Z_2}= A^*_Y\cap\ker(e^{\check\alpha})$, in which case the condition reads $e^{\check\alpha}(\chi^{-1}\delta^{\frac{1}{2}})=1 \Leftrightarrow \chi\in \mathcal Q_{-\alpha}$. 
\item Or otherwise, $^{w_\alpha}(A^*_{Z_1})\ne A^*_{Z_1}$ and we would have $^{w_\alpha}(\Adm_{Z_1})=\Adm_{Z_2}$ if and only if $e^\alpha|_{B_z}=1$, which is impossible.
\end{enumerate}
\end{remark}

\subsection{Convergence and rationality} 

\subsubsection{}\label{sssconvergence} Having examined the question of $(B,\chi^{-1}\delta^{-\frac{1}{2}})$-equivariant distributions on a single $\BB$-orbit, we now examine whether we can extend them to the whole $X$ or, in other words, whether a sequence of the form (\ref{equivdistr}) is surjective on the right. The idea is to use the integral expression (\ref{SY1}), if it converges for all $\phi\in C_c^\infty(X)$, to define an equivariant extension of the distribution $S_\chi^Y$ to the whole space $C_c^\infty(X)$. Then one shows that the resulting morphism is rational in $\chi$, and thus can be extended to almost every $\chi\in \Adm_Y$. As a corollary, for all $\chi$ which do not belong to the ``poles'' of the intertwining operators, we deduce that the sequence (\ref{equivdistr}) is surjective on the right. 

In order to understand the asymptotic behavior of our distributions in the closure of an orbit $Y$, we make use of resolution of singularities. The idea of using resolution of singularities to establish meromorphic properties of certain distributions originates in Atiyah \cite{At} and Bernstein-Gelfand \cite{BG} who used it in the archimedean case. In the $p$-adic case it has been developed and used by Igusa (see \cite{Ig}); for theorems close in formulation to what we need see Denef \cite[Theorem 3.1]{Den} and Deshommes \cite[Th\'eor\`eme 2.5.1]{Dh}. 

Recall that, by Hironaka's embedded resolution of singularities  \cite{Hi}, given a $k$-rational $\BB$-orbit $\YY$ with closure $\overline{\YY}$ there exists a (canonical) regular $k$-scheme $\widetilde{\YY}$ and proper $k$-morphism $p:\widetilde{\YY}\to\overline{\YY}$, which is an isomorphism on $\YY$ and such that the inverse image $\EE$ of $\overline\YY\smallsetminus\YY$ is an effective divisor (the ``exceptional divisor''); moreover, for every point $y_0\in \widetilde\YY(k)$ the geometrically irreducible components of $\EE$ which contain $y_0$ are defined over $k$ and have normal crossings. In other words, the equations of these irreducible divisors around $y_0$ are linearly independent in $\mathfrak m_{y_0}/\mathfrak m_{y_0}^2$, and in particular they form part of a system of coordinates of the Hausdorff topology around $y_0$. Notice that the map $\widetilde Y\to \overline Y$ is surjective on $k$-points, since the proper algebraic morphism induces a proper map in the Hausdorff topology.

In this section it does not make a difference whether $Y$ splits into many $B$-orbits or not, so we will pretend that it does not. To understand the behavior of $S_\chi^Y$, applied to any $\phi\in C_c^\infty(\overline Y)$, it is first better to write it as: $\int_Y \phi(y) d\mu(y)$, where $d\mu$ is a $B$-eigenmeasure on $Y$ with weight $\chi^{-1}\delta^{-\frac{1}{2}}$. In turn, $d\mu$ can be written as $|\omega| |f_1|^{s_1} \cdots |f_m|^{s_m}$ where $\omega$ is a top-degree $\BB$-eigenform and the $f_i$'s are $\BB$-semiinvariants on $\YY$.

We can now pull back $\phi$, the $f_i$'s and $\omega$ to $\widetilde Y$ in order to express the integral as an integral on $\widetilde Y$; the corresponding measure will be denoted by $p^*d\mu$. A divisor $\DD$ on $\widetilde{\YY}$ defines a valuation $v_D$ on rational sections of line bundles, and we can extend the ``exponential'' of this valuation to $d\mu$ or, equivalently, to its weight by setting:
$$q^{v_D(\chi^{-1}\delta^{-\frac{1}{2}})} = q^{v_D(\omega)+s_1 v_D(f_1)+\dots+ s_m v_D(f_m)}.$$
This is a regular function of $\chi\in\Adm_Y$. Notice that the exponent ${v_D(\chi^{-1}\delta^{-\frac{1}{2}})}$ is only well-defined modulo $\frac{2\pi i}{\log q}$.

Let $y_0 \in \widetilde Y\cap \supp(p^*\phi)$. Let $\DD_1, \dots, \DD_k$ be the irreducible components of $\EE$ which contain $y_0$. Then there exist local coordinates $x_1, \dots, x_n$ identifying a neighborhood of $y_0$ in the Hausdorff topology with a neighborhood of $0$ in $k^n$ such that $D_i=\{x_i=0\}$ for $1\le i\le k$. Moreover, $|p^*d\mu|=\prod_{i=1}^{k} |x_i|^{r_i} dx_1\cdots dx_n$ (up to a constant) in this basis, where $r_i=q^{v_{D_i}(\chi^{-1}\delta^{-\frac{1}{2}})}$ . 

Therefore, in a neighborhood of $y_0$, the integral (\ref{SY2}) is equal (up to a constant) to
\begin{equation}\label{local}
\int p^*\phi(x) \prod_{i=1}^{k} |x_i|^{r_i} dx_1 \cdots dx_n
\end{equation}
Recall also that $p^*\phi$ is locally constant. From this we deduce:

\begin{proposition}\label{convergence}
\begin{enumerate}
\item \label{part1}
The integral (\ref{SY2}), representing $S_\chi^Y(\phi)$, converges for all $\phi\in C_c^\infty(\overline{Y})$, for $\chi$ in an open subregion of $\Adm_Y$. It is rational in $\chi$, and its poles are products of factors of the form: 
$$\frac{1}{1-q^{-v_D(\chi^{-1}\delta^{-\frac{1}{2}})-1}}$$ 
where $\DD$ denotes an irreducible component of the exceptional divisor $\EE$ of $\widetilde\YY$. The resulting functional is also to be denoted by $S_\chi^Y$.
\item \label{part2}
Let $\chi\in\Adm_Y$ such that (the rational continuation of) $S_\chi^Y(\phi)$ does not have a pole at $\chi$, let $f\in k[\YY]^{(\BB)}$ of weight $\psi^{-1}$ and let $\phi_i$ denote the restriction of $\phi$ to the set where $|f|=q^{-j}$. For $\kappa\gg 0$, 
$$ S_{\chi\psi^{\kappa}}^Y = \sum_j S_{\chi\psi^{\kappa}}^Y (\phi_j) $$
with the above sum converging absolutely. (Notice that it is not required of the individual summands to be given by a convergent integral.)
\item \label{part3}
Let $A_0$ denote the maximal compact subgroup of $A$ and consider the lattice $A/A_0A_Y$. Choose a point $y\in Y$ and let $t$ denote the map $\phi\mapsto \int_{(A_0U)_y\backslash A_0U} \phi(a \cdot yu) du$ from $C_c^\infty(\overline Y)$ to $C^\infty(A/A_0A_Y)$. Its image is supported on a translate of a cone of the form $\{a |\left<a,\psi\right><\epsilon\}$ for some $\psi\in \Hom(A/A_0A_Y,\mathbb Z), \epsilon>0$ and satisfies
$$ t(\phi)(a) \ll e^{\kappa\left<a,\psi\right>}$$ for some $\kappa$.
\end{enumerate}
\end{proposition}

\begin{remark}
 It is clear that all poles of the form $\frac{1}{1-q^{-v_D(\chi^{-1}\delta^{-\frac{1}{2}})-1}}$, where $\DD$ is an irreducible component of $\EE$ which has a $k$-point, will appear for suitable $\phi$. However, two distinct divisors $\DD_i$ may induce the same $v_{D_i}$, in which case the pole will not necessarily appear with multiplicity two. For instance, in $T\backslash \PGL_2$ ($\TT$ a split torus) there are two colors ($B$-stable prime divisors) which induce the same valuation but do not intersect.
\end{remark}

\begin{proof}
It is obvious that the integral (\ref{local}) converges absolutely if ${v_D(\chi^{-1}\delta^{-\frac{1}{2}})}$ is large enough (in fact, greater than $-1$). Recall that (by the assumption that $\XX$ is quasi-affine) there exists a non-zero $f\in k[\overline{\YY}]^{(\BB)}$ which vanishes on $\overline{\YY} \smallsetminus \YY$. Multiplying by a high enough power of $|f|$ we can achieve the desired valuation for all divisors $D_i$. The rationality and stated form of the poles are immediate from (\ref{local}).

For the second assertion, let $y_0\in\widetilde Y$ be a point as above, and assume that $v_{D_i}>0$ for $i=1,\dots,j$, $v_{D_i}(f)=0$ for $i>j$. Then for $\kappa\gg 0$ the exponents $r_i$ of (\ref{local}) corresponding to $\chi\psi^i$ are strictly increasing affine functions of $\kappa$ for $i\le j$, while for $i>j$ the values of $r_i$ do not depend on $\kappa$. Therefore, if $\phi$ is supported in a neighborhood of $y_0$ then $$S_{\chi\psi^\kappa}^Y(\phi_j) \ll K q^{-\kappa j}\cdot \prod_{i=j+1}^k \frac{1}{1-q^{-r_i-1}}$$
with the constant $K$ depending on $\phi,\chi, f$ but not on $\kappa$ or $j$, and this establishes the claim.

Finally, for the third assertion, notice that the integral under consideration is the integral of $\phi$ restricted to the set where $|f_i|$ have a fixed value, for all $f_i$ in a set of generators for $k[\overline\YY]^{(\BB)}$. Therefore, the estimate follows from the same considerations as above, namely the asymptotic behavior (\ref{local}) of the integral and the existence of an $f\in k[\overline{\YY}]^{(\BB)}$ which vanishes on $\overline{\YY} \smallsetminus \YY$. (Here we use exponential notation for the weight $\psi^{-1}$ of $f$.)
\end{proof}

\subsection{Discussion of the poles} \label{sspoles} 

\subsubsection{} Let $S_\chi^{Y}: C_c^\infty(X)\to I(\chi)$ be as above. By the ``poles'' of $S_\chi^{Y}$ we mean the smallest divisor $M\subset \Adm_Y$ which contains the polar divisors of $S_\chi^{Y}(\phi)$ for all $\phi\in C_c^\infty(X)$. We have proven above that there exists such a divisor, i.e.\ that there is only a finite number of distinct irreducible polar divisors appearing for all $\phi$. 

In Proposition \ref{convergence} we gave a description of the poles of $S_\chi^{Y}$ in terms of geometric data of our spherical variety. There is also a representation-theoretic understanding of the poles, discussed in \cite[2.6]{Ga}, which leads to \emph{necessary} conditions for the poles to appear.

Let $M$ be a closed prime divisor of $A_Y^*$. The local ring $\mathfrak o_{A_Y^*,M}$ is principal. Hence, if $M$ is contained in the subvariety where $S_\chi^{Y}$ has poles, there is $f\in \mathfrak m_{A_Y^*,M}$ (the maximal ideal) such that $S_\chi'^{Y}:= f(\chi) S_\chi'^{Y}$ is \emph{regular} and \emph{nonzero} on a dense subset of $M$.\footnote{In fact, as we saw in Proposition \ref{convergence}, the polar divisors in our case are always principal.} However, the functional $S_\chi'^{Y}$ was regular when restricted to $C_c^\infty(Y)$, so the functional $S_\chi'^{Y}$ will vanish on $C_c^\infty(Y)$ for $\chi\in M$ and will be supported on $\overline Y\smallsetminus Y$. We deduce that, \emph{for $M$ to be a polar divisor of $S_\chi^{Y}$, it has to be contained in the variety of admissible characters of a smaller orbit in the closure of $Y$}. 

Since we already know (by Proposition \ref{convergence} and the remark following it) that some of the possible poles \emph{will} appear, we can extend Garrett's results as follows:

\begin{proposition}
For every irreducible component $\DD$ of $\EE$ (cf.\ \S\ref{sssconvergence}) such that $\DD(k)\ne\emptyset$ there is a $k$-rational $\BB$-orbit $\ZZ$ in  $p(\DD)$ such that \linebreak $\Adm_Z \supset \{\chi| q^{-v_D(\chi^{-1}\delta^{-\frac{1}{2}})}=q\}$.
\end{proposition}

\begin{proposition} \label{pickup}
 Let $\ZZ$ be a $\GG$-orbit on $\XX$. There exists a $k$-rational $\BB$-orbit $\VV$ with $\overline\VV\supset \ZZ$ such that: 
\begin{enumerate}
 \item The inverse image $\widetilde\ZZ:=p^{-1}\ZZ\subset\widetilde\VV$ (in the notation of \S\ref{sssconvergence} with $\YY=\VV$) is an effective divisor and each irreducible component $\DD\subset\widetilde\ZZ$ dominating $\ZZ$ induces the same $v_Z:=v_D\in\Hom(\varchi(\VV),\mathbb Z)$ by restriction to $k(\VV)^{(\BB)}$.
 \item The set $\Adm_Z$ is precisely equal to $\{\chi \in\Adm_V| q^{-v_Z(\chi^{-1}\delta^{-\frac{1}{2}})}=q\}$.
 \item For every $B$-orbit $Z_1\subset Z$ contained in the closure of a $B$-orbit $V_1\subset V$, the rational family of morphisms: $$S_\chi'^{V_1}:=(1-q^{-v_Z(\chi^{-1}\delta^{-\frac{1}{2}})-1})\cdot S_\chi^{V_1} : C_c^\infty(V)\to \CC_{\chi\delta^{\frac{1}{2}}}$$
specializes to $S_\chi^{Z_1}$ for $\chi\in \Adm_Z$.
\end{enumerate}
\end{proposition}

\begin{proof}
We first show that the set of $k$-rational $\BB$-orbits $\VV$ with the property that the inverse image $\widetilde\ZZ:=p^{-1}\ZZ\subset\widetilde\VV$ (in the notation of \S\ref{sssconvergence} with $\YY=\VV$) is an effective divisor is non-empty: We know that $\ZZ$ has a $k$-point, say $z$, in its open $\BB$-orbit. Performing the resolution of singularities as in \S\ref{sssconvergence}, with $\YY=\VV_1:=\mathring\XX$ and $\tilde z_1$ a $k$-point in the preimage of $z$, there is a divisor $\DD_1\subset \widetilde\VV_1$ with $k$-points arbitrarily close to $\tilde z_1$. The image of $\DD_1$ in $\XX$ is contained in an absolutely irreducible, $k$-rational, $\BB$-stable closed subvariety (which is therefore the closure of a $\BB$-orbit $\VV_2$) containing $\ZZ$. We repeat this proccess with $\VV_2$, $\VV_3$ (constructed inductively), etc.\ until $\VV_i=\mathring \ZZ$, which has to occur for reasons of dimension.

Now pick $\VV$ such a $k$-rational $\BB$-orbit of minimal possible dimension.

Now let $Z_1$ denote the $B$-orbit of a point $z\in Z$, and assume that there is a $B$-orbit $V_1\subset V$ containing $Z_1$ in its closure. Let $\tilde z$ be a preimage of $z$ in $\widetilde\VV$ and let $\DD\subset p^{-1}(\widetilde\ZZ)$ be an irreducible $k$-rational divisor containing $\tilde z$ and dominating $\ZZ$. Denote by $v_D$ be the corresponding valuation, considered as an element of $\Hom(\varchi(\VV),\mathbb Z)$ by restriction to $k(\VV)^{(\BB)}$. We have seen that for $\chi$ such that
\begin{equation} \label{chipickedup}
q^{-v_D(\chi^{-1}\delta^{-\frac{1}{2}})}=q
\end{equation}
 the distributions $S_\chi'^{V_1}:=(1-q^{-v_D(\chi^{-1}\delta^{-\frac{1}{2}})-1})\cdot S_\chi^{V_1} : C_c^\infty(V)\to \CC_{\chi\delta^{\frac{1}{2}}}$ specialize to distributions supported on the closure of $Z_1$. We show that in fact they specialize to $S_\chi^{Z_1}$: 

For $\phi\in C_c^\infty(\mathring V \cup \mathring Z)$ and $\chi\in \Adm_Z$ we have $S_\chi'^{V_1} \sim S_\chi^{Z_1}(\phi)$. There is a $\BB$-semiinvariant $f$ on $\XX$ which vanishes precisely on all $\BB$-orbits which do not contain $\ZZ$ in their closure, while it is not identically equal to zero on $\ZZ$ (this is \cite{KnLV}[Corollary 1.7]). By the minimality of $\VV$, the restriction of $f$ to $V$ is zero precisely in the complement of $\mathring V\cup \mathring Z$.
 Applying Proposition \ref{convergence}, (\ref{part2}), with this $f$ and $\psi^{-1}$ denoting the weight of $f$, we have that for large $m$ and $\chi\in \Adm_Z$:
$$ S_{\chi\psi^m}'^{V_1}(\phi) = \sum_j S_{\chi\psi^m}'^{V_1}(\phi_j) \sim \sum_j S_{\chi\psi^m}^{Z_1}(\phi_j) = S_{\chi\psi^m}^{Z_1}(\phi).$$
(Notice that since $\DD$ dominates $\ZZ$, $\chi\psi^m$ also satisfies (\ref{chipickedup}) if $\chi$ does.)
By the rationality of these distributions, we deduce that for every $\chi$ satisfying (\ref{chipickedup}) we have $S_{\chi}'^{V_1} \sim S_{\chi}^{Z_1}$.

It now follows that all $\chi$ satisfying (\ref{chipickedup}) belong to $\Adm_Z$. On the other hand, $\dim\Adm_V-\dim\Adm_Z\ge 1$ (since by Theorem \ref{extension} all $\BB$-semiinvariants on $\ZZ$ extend to $\VV$ and there exist $\BB$-semiinvariants on $\VV$ which vanish on $\ZZ$). Therefore, we deduce that $\Adm_Z=\{q^{-v_D(\chi^{-1}\delta^{-\frac{1}{2}})}=q\}$ and, a posteriori, all such $v_D$ are equal to some $v_Z$. (This is, of course, an unpleasantly indirect proof of this fact.)
\end{proof}

\begin{remark}
The statement is not true for $\ZZ$ any $\BB$-stable set; indeed as we shall see in the next section in the case of $\TT\backslash\PPGL_2$ both closed orbits define the same valuation on $k[\XX]^{(\BB)}$; renormalizing the intertwining operator of $\mathring\XX$ at its pole ``picks up'' a sum of the intertwining operators of the closed orbits.

It is also not true that a pole for $S_\chi^Y$ (where $Y$ is any $B$-orbit) necessarily implies that the corresponding sequence of the form (\ref{equivdistr}) is not surjective on the right. A basic example of this will be encountered in our discussion of Jacquet modules.
\end{remark}

\begin{corollary}
For every $\BB$-orbit $\YY$, we have $\Adm_Y\subset {^w\!\Adm_{\mathring X}}$ for some $w\in [W/W_{P(X)}]$.
\end{corollary}

Recall that $\PP(\XX)$ denotes the standard parabolic $\{g\in\GG| \mathring\XX g =\mathring \XX\}$ and $W_{P(X)}$ the Weyl group of its Levi, and that $[W/W_{P(X)}]$ denotes representatives of minimal length.

\begin{proof}
Assume first that $\YY$ belongs to the open $\GG$-orbit. Let $w_1,w_2,...,w_r$ be simple reflections which succesively raise $\YY$ to $\mathring \XX$, hence $\codim\YY=r$. It is known that $w=w_1 w_2\cdots w_r \in [W/W_{P(X)}]$ \cite[Lemma 5.(iii)]{BrOr}. From the discussion of \S\ref{sssadmissible}, $\Adm_Y\subset {^w\!\Adm_{\mathring X}}$.

Now let $\YY=\mathring \ZZ$, where $\ZZ$ is a smaller $\GG$-orbit. Let $\VV$ be as in Proposition \ref{pickup}. It follows that $\Adm_Y\subset\Adm_V$.

For a general $B$-orbit the claim now follows by applying the above two steps and the fact that $\PP(\XX)\subset\PP(\ZZ)$ for every $\GG$-orbit $\ZZ$.
\end{proof}

The importance of this result will be that most information about the unramified spectrum of a spherical variety can already be retrieved by looking at the open orbit.

Combining all the results above, we have proven the following:

\begin{theorem} \label{existence}
Assume that $\mathring X$ carries a $B$-invariant measure (cf.\ \S\ref{invariantmeasure}). A necessary condition for the existence of a non-zero morphism: $C_c^\infty(X)\to I(\chi)$ is that: 
\begin{equation}\label{condchi}\chi\in {^w\!\left(\delta^{-\frac{1}{2}}A_X^*\right)} \textrm{ for some } w\in [W/W_{P(X)}].\end{equation}
 For every $B$-orbit $Y$ on $X$ and $\chi\in\Adm_Y$ there exists a natural (up to a non-zero regular function on $\Adm_Y$) family $S_\chi^Y: C_c^\infty(X)\to I(\chi)$, rational in $\chi$. For almost all $\chi$ satisfying the condition (\ref{condchi}), the space of morphisms $C_c^\infty(X)\to I(\chi)$ admits a basis consisting of all $S_\chi^Y$ with $\Adm_Y\ni\chi$.
\end{theorem}

Finally, we comment on the ``shift'' $\delta^{-\frac{1}{2}}$ which appears in the description of admissible characters. We have seen a typical example where this shift has a significance (i.e.\ is not absorbed by $A_X^*$): In the variety $\PGL_2\backslash\PGL_2$ whose spectrum consists only of the trivial representation (which is a subrepresentation of $I(\delta^{-\frac{1}{2}})$). In fact, this is essentially the only appearance of a non-trivial shift:

\begin{lemma}
 Under the assumptions of the above theorem, the variety $\delta^{-\frac{1}{2}} A_X^*$ is equal to $\delta_{L(X)}^{-\frac{1}{2}} A_X^*$.
\end{lemma}

Here $\LL(\XX)$ denotes the Levi of $\PP(\XX)$ and $\delta_{L(X)}$ is the modular character of its Borel subgroup, in other words, it is equal to $e^{2\rho_{L(X)}}$ where $\rho_{L(X)}$ is the half-sum of positive roots of $\LL(\XX)$. Notice that $\rho_{L(X)}$ is orthogonal to $\mathfrak a_X^*$ by Lemma \ref{charlemma}.

\begin{proof}
 We are using the following two facts: First, by assumption $\mathring \XX$ carries a $\BB$-invariant measure. This implies that $\delta|_{A_X}= \delta_{B_x}$ where $x$ is a point on $\mathring X$. Secondly, by \cite{Po}, the open $\BB$-orbit of a spherical variety $\XX$ is $\BB$-isomorphic to the open $\BB$-orbit of a horospherical variety $\SS$ (one whose stabilizer contains a maximal unipotent subgroup) with $\PP(\SS)=\PP(\XX)$. It follows that $\delta_{B_x}=\delta_{L(X)}$.
\end{proof}

\subsection{Jacquet modules}

Any $B$-equivariant functional $V\to \CC_{\chi\delta^{\frac{1}{2}}}$ (for $(\pi,V)$ a smooth representation of $G$) factors through the \emph{Jacquet module} $V_U$: This is, by definition the maximal quotient of $V$ where $U$ acts trivially; equivalently, it is equal to the quotient of $V$ by the span of $\{v - \pi(u)v| u\in U, v\in V\}$. It is well-known that the $A$-equivariant functor $V\mapsto V_U$ is exact, due to the fact that $U$ is filtered by compact subgroups. In fact, since we are only considering unramified principal series, we may as well compose with the functor $V_U\mapsto V_{A_0U}$ ((co-)invariants for the maximal compact subgroup $A_0$ of $A$), which is also exact; we will call $V_{A_0U}$ the \emph{unramified} Jacquet module.

In what follows we examine the Jacquet modules for some basic $\GL_1$-and $\GL_2$-spherical varieties. We do this in order to demonstrate how the method of intertwining operators gives us information on the Jacquet module; to show that the Jacquet module does not, in general, have a very simple geometric description; and to discuss what happens at characters $\chi$ on the poles of the intertwining operators, where the above method fails to prove surjectivity of (\ref{equivdistr}) on the right.

\begin{example}\label{ex1} Let $\XX=\mathbb A^1$, as a $\GGL_1$-spherical variety. From the two orbits $\mathring X=k^\times$ and $Z=\{0\}$ we have the sequence:
$$ 0\to C_c^\infty(k^\times) \to C_c^\infty(k) \to \CC \to 0.$$
The corresponding sequence of unramified Jacquet modules is:
\begin{equation}\label{A1}
 0\to C_c(\mathbb Z) \to C_c^\infty(\mathbb Z\cup\{-\infty\}) \to \CC \to 0
\end{equation}
where $C_c(\mathbb Z)$ denotes compactly supported sequences on $\mathbb Z$ and $C_c^\infty(\mathbb Z\cup\{-\infty\})$ denotes sequences supported away from $-\infty$ which stabilize in a neighborhood of $-\infty$.

Tate's thesis shows that the intertwining operator $S_\chi^{\mathring X}$ has a pole at $\chi=1$, which is exactly where the intertwining operator $S_\chi^Z$ appears. (This is a general phenomenon which will be discussed below.) The element $1-\chi$ of the Bernstein center (respectively: the Hecke algebra element which maps the sequence $(a_n)_n$ to the sequence $(a_n-a_{n-1})_n$) is clearly injective from $C_c^\infty(k)$ onto $C_c^\infty(k^\times)$ (resp.\ on the corresponding unramified parts), and this implies that the Jacquet module of $C_c^\infty(k)$ is isomorphic to that of $C_c^\infty(k^\times)$. Its unramified part is, as mentioned, $C_c(\mathbb Z)\simeq \CC[T,T^{-1}]$ and if we rewrite (\ref{A1}):
$$ 0\to \CC[T,T^{-1}]\to\CC[T,T^{-1}]\to\CC\to 0$$
then the map on the right is just evaluation at $T=1$. 

We deduce, in particular, that the sequence (\ref{equivdistr}) is not surjective on the right in this case.
\end{example}

\begin{example}\label{ex2}
Let $\XX=\TT\backslash \PPGL_2$, where $\TT=\AA$ is a $k$-split torus. Let $\YY_1,\YY_w$ denote the two closed $\BB$-orbits represented by the elements $1$ and $w$ (a representative for the non-trivial Weyl group element), respectively. By the Bruhat decomposition, $\PPGL_2 = \BB \sqcup \BB w\BB$, the orbit $\YY_w$ has an open, $\BB$-stable neighborhood $\AA\backslash(\BB w\BB) = \AA\backslash(\BB w \UU) \simeq \mathbb A^1 \times \mathbb A^1$ where the action of $\BB=\AA\UU$ is described as follows: $\AA$ acts by the character $\mathfrak d$ on the first factor and by $\mathfrak d^{-1}$ on the second, and $\UU\simeq\Ga$ acts by translations on the second. The geometric quotient of this open set by the action of $\UU$ exists and is equal to $\mathbb A^1$ by projection onto the first factor. It is easy to see that integration over the orbits of $U$ defines an isomorphism: $C_c^\infty(A\backslash BwU)_U \simeq C_c^\infty(k)\otimes \CC_{\delta}$, where the action of $A$ on $k$ is via the character $\delta$. There is a $\GG$-automorphism (multiplication on the left by $w$, the non-trivial element of the Weyl group) which carries one closed orbit to the other and the open neighborhood $\AA\backslash \BB w\UU$ to the open neighborhood $\AA\backslash w\BB w\UU$. Therefore, we have:
$$C_c^\infty(X)_U\simeq (C_c^\infty(k)\oplus C_c^\infty(k))/C_c^\infty(k^\times)^{\operatorname{diag}}\otimes \CC_\delta$$
where the diagonal copy of $C_c^\infty(k^\times)\otimes \CC_\delta$ corresponds to the $U$-coinvariants of $C_c^\infty(\mathring X)$ and is embedded with a minus sign in one of the factors. In other words, up to twisting by $\CC_\delta$ the Jacquet module corresponds to functions $f$ on the non-separated ``affine line with doubled origin'' $k^\times \cup\{0_1,0_2\}$ which are locally constant on $k^\times$, vanish eventually as $x\to\infty$, while as $x\to 0$ they eventually stabilize to $f(0_1)+f(0_2)$.

Using the previous example, it follows that the sequence (\ref{equivdistr}) is not surjective on the right; at $\chi=\delta^{\frac{1}{2}}$ the dimension of intertwining operators is equal to two, coming from the two closed orbits.
\end{example}

\begin{example}\label{ex3}
Let now $\XX=\UU\backslash\PPGL_2$. The Bruhat decomposition gives us a filtration:
$$0\to C_c^\infty(U\backslash BwB) \to C_c^\infty(X) \to C_c^\infty(U\backslash B)\to 0$$
with corresponding Jacquet modules (cf.\ \cite[Proposition 6.2.1]{Cas}):
$$0\to C_c^\infty(k^\times)\to C_c^\infty (X)_U \to C_c^\infty (k^\times) \to 0$$
and unramified Jacquet modules:
$$0\to \CC[T,T^{-1}] \to C_c^\infty(X)_{A_0U} \to \CC[T,T^{-1}]\to 0.$$
The latter is a sequence of $\mathcal H(G,K)=\CC[T,T^{-1}]$-modules (where convolution with elements of the Hecke algebra corresponds to multiplication of polynomials), and since this ring is a principal ideal domain and the modules are free, the sequence splits, so we have (non-canonically):
$$ C_c^\infty(X)_{A_0U}\simeq \CC[T,T^{-1}] \oplus \CC[T,T^{-1}].$$
Therefore we see that although the intertwining operator $S_\chi^{\mathring X}$ has a pole at $\chi=1$, the corresponding sequence (\ref{equivdistr}) is surjective on the right in this case, and the dimension of intertwining operators is constantly equal to two.

\end{example}

In any case, our intertwining operators are enough to characterize the image of a $\phi$ in the Jacquet module:

\begin{lemma}\label{enough}
A vector $\phi\in C_c^\infty(X)$ lies in the kernel of the Jacquet morphism $C_c^\infty(X)\to C_c^\infty(X)_U$ if and only if the integral of $\phi$ over every horocycle ($U$-orbit) is zero. Similarly, it lies in the kernel of the unramified Jacquet morphism if and only if its integral over every $A_0U$-orbit is zero. The integral of $\phi$ over all $A_0U$-orbits in a given $B$-orbit $Y$ is zero if and only if $S_\chi^Y(\phi)=0$ for generic $\chi\in\Adm_Y$.
\end{lemma}

\begin{proof}
If $Y$ is a $B$-orbit and $\phi \in C_c^\infty(Y)$ is compactly supported, then it is obvious that its image in $C_c^\infty(Y)_U$ (resp.\ $C_c^\infty(Y)_{A_0U}$) is zero if and only if its integral over all orbits of $U$ (resp.\ $A_0U$) is zero.

Now let $\phi\in C_c^\infty(X)$. Since the kernel of the Jacquet morphism is generated by elements of the form $f-R(u) f$, where $R$ denotes the right regular representation, it is clear that if $\phi$ belongs to the kernel then its integral over every horocycle is zero. Conversely, using the standard filtrations of Jacquet modules and by induction on the orbit dimension, we will prove that $\phi$ lies in the kernel of the Jacquet morphism if all its integrals on horocycles are zero. So, assume that the integral of $\phi$ over every horocycle is zero and let $m$ be the minimal dimension of an orbit which intersects the support of $\phi$. Then the image of $\phi$ under
$$ C_c^\infty(X)_U \twoheadrightarrow C_c^\infty(\cup_{\dim Y\le m}Y)_U $$
is zero. The kernel of the above map is equal to the Jacquet module of \linebreak $C_c^\infty(\cup_{\dim Y> m} Y)$, hence $\phi$ differs from a $\phi'\in C_c^\infty(\cup_{\dim Y>m}Y)$ by a function of the form $f-R(u) f$. Since the latter has integral zero over any $U$-orbit, we reduce the problem to $\phi'$, which allows us to complete the proof by induction. The claim about the unramified Jacquet module follows similarly. 

For the last claim, the direction $\Rightarrow$ is, again, obvious. For the inverse, use part (\ref{part3}) of Proposition \ref{convergence}: Multiplying $t(\phi)$ by a suitable character of $B$, it lands in $L^1(A/A_0A_Y)$. By standard Fourier analysis on the discrete abelian group $A/A_0A_Y$, if its Fourier transform is zero then the function itself is zero.
\end{proof}

The importance of the above lemma is that in order to establish certain results we do not have to worry about intertwining operators which may not be expressible in terms of our $S_\chi^Y$'s.

\subsection{Non-trivial line bundles and standard intertwining operators}\label{ssbundle}

As mentioned above, exactly the same arguments apply to intertwining operators: $C_c^\infty(X,\mathcal L_\psi)\to I(\chi)$, where $\psi$ is some complex character of $H$. The condition of admissibility with respect to a $\BB$-orbit $\YY$ is now: $\chi^{-1}\delta^{\frac{1}{2}}|_{B_y}=\delta_{B_y}{^y\psi}^{-1}|_{B_y}$, where $^y\psi$ denotes the character by which the stabilizer of $y$ (a conjugate of $H$) acts on the fiber of the map $\GG\to\HH\backslash\GG$ over $y$.

As a special case of this, the filtration of $\BB\backslash\GG$ defined by the Bruhat decomposition gives rise to the standard intertwining operators for unramified principal series:
$$T_w : I(\chi) \to I({^w\chi})$$
which are rational in $\chi\in A^*$ and are given by the rational continuation of the integral:
$$ \int_{\prod_{\alpha>0,w^{-1}\alpha<0}U_\alpha} \phi(w^{-1}u) du. $$
The above integral expression depends on the choice of a representative for $w$ in $\mathcal N(\AA)$, but only up to a character of $A^*$, therefore we will ignore this dependence whenever we can. Their poles are a union of ``irregular'' divisors as described in \S \ref{ssups}, and one can verify that those are the characters where a smaller Schubert cell can support an intertwining operator into $I(\chi)$. Note that in the case of the variety $\BB\backslash\GG$ Knop's action translates to the action of $W$ on itself by left multiplication.


\section{Interpretation of Knop's action}\label{secKnop}

\subsection{Avoidance of ``bad'' divisors}

The object of this section is to investigate what happens when one composes the morphisms $S_\chi^{Y,*}$ with the intertwining operators for principal series $T_w$. Before we do that, we need to examine issues that might arise from the set of our characters $\chi$ being contained in some of the hypersurfaces where $I(\chi)$ is reducible (and where some of the $T_w$ annihilate a subrepresentation). It turns out that this can only happen for trivial reasons. These trivial reasons are best exhibited in the example of the $\SSL_2$-spherical variety of type G, namely $X$=a point. Then $C_c^\infty(X)$ is the trivial representation of $G$; it is contained as a proper subrepresentation in $I(\chi)$ where $\chi=\delta^{-\frac{1}{2}}$ 
and as a quotient (but not a subrepresentation) in $I({^w\chi})$. Hence, $S_\chi$ maps into that subrepresentation of $I(\chi)$ and $T_w: I(\chi) \to I({^w\chi})$ annihilates its image. As it turns out, this is essentially the only way things can go wrong. We recall from \S \ref{ssups} the definition of the ``bad'' divisors $R_{\check\alpha}=\{\chi\in A^* | \chi^{\check\alpha}=1\}$ and $Q_{\check\alpha}=\{ \chi\in A^* | \chi^{\check\alpha}=q\}$, where $\check\alpha$ is a coroot and $q$ is the order of the residue field.

\begin{lemma}\label{regular}
The variety $\delta^{-\frac{1}{2}} A_X^*$ is never contained in one of the ``irregular'' divisors $R_{\check\alpha}$. It is contained in $Q_{\check\alpha}$ if and only if $\alpha$ is a simple, positive root of the Levi of $\PP(\XX)$.
\end{lemma}

\begin{proof}
Since $\delta^{-\frac{1}{2}}\in \delta^{-\frac{1}{2}}A_X^*$ and $\delta^{-\frac{1}{2}}$ is regular, $\delta^{-\frac{1}{2}}A_X*$ is not contained in any of the $R_{\check\alpha}$.

We have: $\delta^{-\frac{1}{2}} A_X^*\subset Q_{\check\alpha} \Rightarrow e^{-\rho}\in Q_{\check\alpha} \iff \left<\rho,\check\alpha\right>=1 \iff \alpha\in\Delta$. In that case, we see that $w_\alpha$ has to centralize $A_X^*$ which, by non-degeneracy (\S \ref{ssnondeg}), implies that $\alpha\in \Delta_{P(X)}$. The converse is easily checked.
\end{proof}

\begin{corollary}\label{irreducible}
For generic $\chi\in \delta^{-\frac{1}{2}}A_X^*$ the image of $S_\chi^{\mathring X}$ in $I(\chi)$ is irreducible.
\end{corollary}

\begin{proof}
Indeed, since a generic $\chi$ is contained only in those $Q_{\check\alpha}$ with $\alpha$ simple, positive and appearing in the Levi quotient $\LL(\XX)$ of $\PP(\XX)$, and since the stabilizer inside $\PP(\XX)$ of a generic point contains, modulo the unipotent radical of $\PP(\XX)$, the commutator subgroup of $\LL(\XX)$, it follows that for such $\chi$ the image of $S_\chi^{\mathring X}$ in $I(\chi)$ belongs to the irreducible subspace induced from the trivial representation of the commutator of $\LL(\XX)$.
\end{proof}

\subsection{The basic theorem}

\begin{theorem}\label{maintool}
 Let $\YY$ be a $k$-rational $\BB$-orbit on $\XX$ and $\alpha$ a simple root such that $\YY$ is of maximal rank in $\YY\cdot\PP_\alpha$. The following describes the composition of $T_{w_\alpha}$ with elements of $\underline S_\chi^Y$ (cf.\ \S\ref{ssweighted}), for $\chi\in\Adm_Y$:
\begin{enumerate}
\item If $(\YY,\alpha)$ is of type $G$ then $T_{w_\alpha}\circ S_1 = 0$ for every $S_1\in \underline S_\chi^Y$.
\item If $(\YY,\alpha)$ is of type $U$ or $T$ then $T_{w_\alpha}\circ S_{\tilde\chi}^Y \sim S_{^{w_\alpha}{\tilde\chi}}^{^{w_\alpha}Y}$.
\item If $(\YY,\alpha)$ is of type $N$ then $T_{w_\alpha}\circ S_{\tilde\chi,\zeta}^Y \sim S_{^{w_\alpha}\tilde\chi,\zeta}^Y$.
\end{enumerate}
Moreover, in cases $T$ and $N$, if $\ZZ$ is a smaller rational orbit in $\YY\PP_\alpha$ then for generic $\chi\in \Adm_Z$ we have $T_{w_\alpha}\circ S_{\tilde\chi}^Z\sim S_{^{w_\alpha}\tilde\chi}^Y$ (resp.\ $T_{w_\alpha}\circ S_{\tilde\chi}^Z\sim S_{^{w_\alpha}\tilde\chi,\zeta}^Y$ where $\zeta$ is the coset corresponding to a split torus in case $N$).
\end{theorem}

The proof will be performed in two steps: First we will show it for functions whose support on $(\YY\cdot\PP_\alpha)(k)$ is compact, and then we will extend it to all $\phi\in C_c^\infty(X)$.

\begin{proposition} \label{step1}
 The statement of Theorem \ref{maintool} is true when the $S_*^*$'s, $T_{w_\alpha}\circ S_*^*$'s are viewed as functionals and restricted to $\phi\in C_c^\infty(X)$ with $\phi|_{(\YY\PP_\alpha)(k)}\in C_c^\infty((\YY\PP_\alpha)(k))$.
\end{proposition}

\begin{proof}
If $(\YY,\alpha)$ is of type $G$ then the image of every $S_1\in\underline S_\chi^Y$ is contained in $\Ind_{P_\alpha}^G(\chi\delta^{\frac{1}{2}})\subset I(\chi)$, and $T_{w_\alpha}$ annihilates that subspace.

In each of the other cases, let $S_1\in \underline S_\chi^Y$ be one of the basis elements as in the statement of Theorem \ref{maintool}, according to the type of $(\YY,\alpha)$. 
Let $y\in Y$ and, in the case of type $N$, $y\in\zeta$ where $S_1=S_{\tilde\chi,\zeta}$. Let $\HH_\alpha=(\GG_y\cap \PP_\alpha \mod \UU_{P_\alpha})\subset\LL_\alpha$.
We can write the functional $S_1$ as (the rational continuation of):
$$ S_1 (\phi) = \int_{U_y\backslash U\times A''} \phi(y\cdot u a) \tilde\chi^{-1}\delta^{\frac{1}{2}}(a) du da$$
where $A''=A'^Y$ in cases $T,U$ and $A''=A'^Y_\alpha$ in case $N$. (Notation and abuse of notation as in \S\ref{ssweighted}.)

We first want to reduce the case of type $N$ to the case of type $T$, and the basis $S_{\tilde\chi,\zeta}^Y$ was chosen precisely for that purpose. Namely, consider the quotient $\LL_\alpha\to \LL_\alpha/\ZZ_\alpha\simeq\PPGL_2$; the image of $\HH_\alpha$ is equal to $\HH_2:=\mathcal N(\TT_2)$ for some torus $\TT_2$. Let $\HH_\alpha^1$ denote the preimage of the connected component of $\HH_2$, namely the preimage of $\TT_2$. Then we have a quotient map: $\HH_\alpha^1\backslash\LL_\alpha\to\HH_\alpha\backslash\LL_\alpha$ and the image of $(\HH_\alpha^1\backslash\LL_\alpha)(k)$ intersected with $Y$ corresponds exactly to the coset $\zeta$.

In cases $U$ and $T$ let $\HH_\alpha^1=\HH_\alpha$. The central observation is that in all cases $\HH_\alpha^1\backslash \LL_\alpha$ is a homogeneous spherical variety for a group $\widetilde{\LL_\alpha}$ which acts transitively on its $k$-points, and whose Borel subgroup $\widetilde{B_\alpha}$ acts transitively on the $k$-points of the open orbit. Indeed, in cases $T$ and $N$ we have $\HH_\alpha\cap\BB\subset\ZZ_\alpha$ so we can let $\widetilde{\LL_\alpha}$ be the group $\LL_\alpha/(\HH_\alpha\cap\BB)$. In case $U$, $\AA\cap\HH_\alpha$ normalizes $\HH_\alpha$ so we can let $\widetilde{\LL_\alpha}$ be the group $(\AA/\AA\cap\HH_\alpha)\times\LL_\alpha$ (with the first factor acting ``on the left''). Now it is clear that $\widetilde{B_\alpha}$ acts transitively on the $k$-points of the open orbit, and by Lemma \ref{nbhd} so does $\widetilde{L_\alpha}$ on the $k$-points of $\HH_\alpha^1\backslash \LL_\alpha$. Let $\widetilde{\HH_\alpha^1}$ be the corresponding isotropy group. Notice that $\widetilde{\HH_\alpha^1}\backslash\widetilde{\LL_\alpha}$ is of type $T$ or $U$. Moreover, $\tilde\chi$ can be considered as a character (possibly ramified) of $\widetilde{B_\alpha}$.

Now $S_1$ can be considered as a morphism: $C_c^\infty((\YY\PP_\alpha)(k))\to \Ind_{\widetilde{B_\alpha}}^{\widetilde{L_\alpha}}(\tilde\chi\delta^\frac{1}{2})$ and analyzed into the composition of two morphisms: First, integration over ${U_{P_\alpha}}_y\backslash U_{P_\alpha}$:
\begin{equation}\label{iota}
 f\mapsto \int_{{U_{P_\alpha}}_y\backslash U_{P_\alpha}} f(y\cdot u \cdot \bullet) du
\end{equation}
defines a morphism:
$$ \iota_\alpha: C_c^\infty((\YY\PP_\alpha)(k))\to C_c^\infty(\widetilde{H_\alpha^1}\backslash\widetilde{L_\alpha},\mathcal L_{\delta_{P_\alpha}\delta_{P_y\cap U_{P_\alpha}}^{-1}}).$$

This is followed by integration over $\bar y\cdot \widetilde{B_\alpha}$, where $\bar y \in \widetilde{H_\alpha^1}\backslash\widetilde{L_\alpha}$ is a point mapping to $(y\mod U_{P_\alpha})$:
$$\widetilde{S_{\tilde\chi}^Y}: C_c^\infty(\widetilde{H_\alpha^1}\backslash\widetilde{L_\alpha},\mathcal L_{\delta_{P_\alpha}\delta_{P_y\cap U_{P_\alpha}}^{-1}}) \to \Ind_{\widetilde{B_\alpha}}^{\widetilde{L_\alpha}}(\tilde\chi\delta^\frac{1}{2}).$$

We analyze the composition of $\widetilde{S_{\tilde\chi}^Y}$ with $T_{w_\alpha}$. First, we notice that $T_{w_\alpha}\circ \widetilde{S_{\tilde\chi}^Y}$ has image in $\Ind_{\widetilde{B_\alpha}}^{\widetilde{L_\alpha}}(^{w_\alpha}\tilde\chi\delta^\frac{1}{2})$.

\subsubsection{Cases T and N:} Here, by the fact that $\widetilde B_\alpha$ has a unique open orbit, it follows that for generic $\chi\in\Adm_Y$ there is a unique morphism:
$$C_c^\infty(\widetilde{H_\alpha^1}\backslash\widetilde{L_\alpha},\mathcal L_{\delta_{P_\alpha}\delta_{P_y\cap U_{P_\alpha}}^{-1}}) \to \Ind_{\widetilde{B_\alpha}}^{\widetilde{L_\alpha}}(^{w_\alpha}\tilde\chi\delta^\frac{1}{2}).$$
It follows that $T_{w_\alpha}\circ \widetilde{S_{\tilde\chi}^Y} \sim \widetilde{S_{^{w_\alpha}\tilde\chi}^Y}$ and hence $T_{w_\alpha}\circ S_{\tilde\chi}^Y \sim S_{^w{\tilde\chi}}^{^wY}$ (in case T) and $T_{w_\alpha}\circ S_{\tilde\chi,\zeta}^Y \sim S_{\tilde\chi,\zeta}^Y$ (in case N).

The statement about $T_{w_\alpha} \circ S_{\tilde\chi}^Z$ follows from the remark in \S\ref{sssadmissible}: Since ${^{w_\alpha}\Adm_Z}$ is not contained in $\Adm_Z'$ for any non-open orbit $\ZZ\subset \YY\PP_\alpha$, for generic $\chi\in\Adm_Z$ we must have $T_{w_\alpha}\circ S_{\tilde\chi}^Z\sim S_{^{w_\alpha}\tilde\chi}^Y$. Notice that $T_{w_\alpha}$ is a well-defined and non-zero for generic $\chi\in \Adm_Z$.

\subsubsection{Case U:} Without loss of generality, since $T_{w_\alpha}\circ T_{w_\alpha} \sim \operatorname{Id}$, let $\YY$ be the closed orbit in $\YY\PP_\alpha$. 
Then $\widetilde{S_{\tilde\chi}^Y}$ is given by the functional: 
$$\phi\mapsto \int_{\widetilde{A_Y}\backslash\widetilde A} \phi(\bar y \cdot a) \tilde\chi^{-1}\delta^{-\frac{1}{2}}(a) da.$$
This converges absolutely for every $\tilde\chi$ if $\phi\in C_c^\infty(\widetilde{H_\alpha^1}\backslash\widetilde{L_\alpha},\mathcal L_{\delta_{P_\alpha}\delta_{P_y\cap U_{P_\alpha}}^{-1}}) $, and in the domain of convergence of $T_{w_\alpha}$ we get:
$$T_{w_\alpha}\circ \widetilde{S_{\tilde\chi}^Y}(\phi) =  \int_{U_2}\int_{\widetilde{A_Y}\backslash\widetilde A} \phi(\bar y \cdot awu) \tilde\chi^{-1}\delta^{-\frac{1}{2}}(a) da du$$
which is precisely equal to $\widetilde{S_{^{w_\alpha}\tilde\chi}^{^{w_\alpha}Y}}$.

This completes the proof of the proposition.
\end{proof}

\subsubsection{Orbits in the closure do not contribute.} \label{ssoutside}
To conclude the proof of Theorem \ref{maintool}, we need to show that what we just proved for $\phi$ compactly supported on $(\YY\cdot\PP_\alpha)(k)$ continues to hold for $\phi$ supported in its closure. We now use $S_\chi^Y$ to denote any of the basis elements in the formulation of the theorem, according to the type of $(\YY,\alpha)$ that we are considering. The idea is to split the intersection of the support of $\phi$ with $P_\alpha$ into infinitely many compact pieces, let $\phi_i$ denote the restriction of $\phi$ to the $i$-th piece by $\phi_i$ (hence $\phi=\sum_i \phi_i$ when restricted to $(\YY\cdot\PP_\alpha)(k)$) and use the fact that $S_\chi^{Y}(\phi)=\sum_i S_\chi^{Y}(\phi_i)$ when $\chi$ is such that the integral expression for $S_\chi^{Y}$ converges. The problem is that $T_{w_\alpha}$ and $S_\chi^Y$ will not, in general, converge simultaneously so we cannot use their integral expressions to prove directly that $T_{w_\alpha} \sum_i S_\chi^{Y}(\phi_i) = \sum_i T_{w_\alpha} S_\chi^{Y}(\phi_i)$. To solve this problem, we could make use of the asymptotic estimates of Proposition \ref{convergence}, part (\ref{part2}), with a suitable $f$ (as in Proposition \ref{pickup}). However, asymptotic estimates are unnecessary here:

\begin{lemma}
Let $K_1$ be an open compact subgroup of $P_\alpha$. Let $g_1,\dots,g_m$ be representatives for the orbits of $K_1$ on $B\backslash P_\alpha$. Then there are rational functions $r_1,\dots r_m$ of $\chi$ such that for $\phi\in I(\chi)^{K_1}$ we have: $$ T_{w_\alpha} \phi = \sum_{j=1}^m r_j(\chi) \phi(g_j).$$
\end{lemma}

The lemma is a direct consequence of the rationality of $T_w$ and the fact that $\Ind_B^{P_\alpha} (\chi\delta^\frac{1}{2})^{K_1}$ is finite-dimensional. It is important that we only fix a compact open subgroup of $P_\alpha$, not of the whole group $G$.

Now, given $\phi\in C_c^\infty(\overline{(\YY\cdot\PP_\alpha)}(k))$ fix a compact-open $K_1\subset P_\alpha$ such that $\phi$ is $K_1$-invariant and representatives $g_1,\dots,g_m$ as above and enumerate the $K_1$-orbits on $(\YY\cdot\PP_\alpha)(k)$: $O_1, O_2, O_3, \dots$. Let $\phi_i= \phi \cdot \mathbf{1}_{O_i} \in C_c^\infty((\YY\cdot\PP_\alpha)(k))$. 

Notice that for $g\in P_\alpha$ the sets $O_i g$ define a partition of $(\YY\cdot\PP_\alpha)(k)$ in $g^{-1}K_1 g$-orbits.

For $\chi$ in the region of convergence of the integral expression for $S_\chi^Y$ we have: $S_\chi^Y(R(g_j)\phi)=\sum_i S_\chi^Y(R(g_j)\phi_i)$ for every $j$. Using the previous lemma:
\begin{eqnarray*}
T_{w_\alpha} \sum_i S_\chi^Y(\phi_i) = \sum_{j=1}^m r_j(\chi) \sum_i R(g_j) S_\chi^Y(\phi_i) =\\ = \sum_i \sum_{j=1}^m r_j(\chi) R(g_j) S_\chi^Y(\phi_i) = \sum_i T_{w_\alpha} S_\chi^Y(\phi_i).
\end{eqnarray*}

Using Proposition \ref{step1}, we have $T_{w_\alpha}\circ S_\chi^{Y} (\phi_i)=S_{^w \chi}^{{^w Y}}(\phi_i)$. Hence: $T_{w_\alpha} S_\chi^Y = S_{^w\chi}^{{^w Y}}$. This completes the proof of Theorem \ref{maintool}.

\subsection{Corollaries and examples}

We discuss the implications of Theorem \ref{maintool} for elements of the Weyl group of length greater than one:

\begin{corollary}\label{nezero}\label{cor1}
$T_w\circ \underline S_\chi^{\mathring X} \ne 0$ if and only if $w\in [W/W_{P(X)}]$. For $w\in [W/W_{P(X)}]$, $T_w\circ \underline S_\chi^{\mathring X} = \underline S_{^w\chi}^{{^w\mathring X}}$. If there are no orbits $\YY$ of maximal rank and simple roots $\alpha$ such that $(\YY,\alpha)$ is of type $N$, then  $T_w\circ  S_{\tilde\chi}^{\mathring X} \sim \underline S_{^w\tilde\chi}^{{^w\mathring X}}$ for every $w\in [W/W_{P(X)}]$.
\end{corollary}

\begin{proof}
Every $w=w_1 w_2$ with $w_1\in [W/W_P]$ and $w_2\in W_{P(X)}$ (uniquely). It follows from Theorem \ref{maintool} that $T_{w_2}\circ S_\chi^{\mathring X} = 0$. The elements of $[W/W_{P(X)}]$ are characterized by the fact that $w\alpha>0$ for every (simple) positive root $\alpha$ in the Levi of $\PP(\XX)$. From Lemma \ref{regular} and the properties of intertwining operators (\S \ref{ssups}) it follows that $T_{w_1}$ is an isomorphism for almost every $\chi$ on $\delta^{-\frac{1}{2}} A_X^*$. The second statement follows immediately from Theorem \ref{maintool}.
\end{proof}

\begin{remark}
For simply-laced groups, Brion \cite{BrOr} has shown that the condition ``there are no orbits $\YY$ of maximal rank and simple roots $\alpha$ such that $(\YY,\alpha)$ is of type $N$'' is equivalent to the condition ``there is no simple root $\alpha$ such that $(\mathring\XX,\alpha)$ is of type $N$''.
\end{remark}

By \cite{KnOrbits}, the stabilizer of the open orbit is $W_X\ltimes W_{P(X)}$. Moreover, by definition $W_X\subset [W/W_{P(X)}]$. The points of $A_X^*$ are left stable by $W_{P(X)}$. Hence $T_w \circ \underline S_\chi^{\mathring X} = \underline S_{^w\chi}^{\mathring X}$ (for generic $X$) if and only if $w \in W_X$.

Denote by $S_{\tilde\chi}$ the operator $S_{\tilde\chi}^{\mathring X}$ for $\chi\in\delta^{-\frac{1}{2}}A_X^*$. The problem that if there are $(\YY,\alpha)$ of type $N$ then we cannot explicitly ``diagonalize'' the composition of elements of $\underline S_\chi:= \underline S_\chi^{\mathring X}$ with $T_w$, for $w$ of length greater than one, can be amended non-explicitly as follows: We claim that there still exists a rational basis $(S^i_\chi)_i$, $i=1,\dots,|H^1(k,\AA_X)|$, of $\underline S_\chi$ such that $T_w\sim S^i_\chi \sim S^i_{^w\chi}$ for every $i$ and $w\in W_X$. Indeed, we may re-normalize the operators $T_w$ so that they satisfy $T_{w_1}\circ T_{w_2} = T_{w_1 w_2}$ (for instance, as equivariant Fourier transforms on $U\backslash G$, cf.\ \cite{BK}). Then the matrices $b_w(\chi)$ of the relation: $T_w(\chi) [S_{\tilde\chi}]_{\tilde\chi} = b_w(\chi) [S_{\tilde\chi}]_{\tilde\chi}$ are 1-cocycles from $W_X$ to $\GL_n(\CC(\delta^{-\frac{1}{2}}A_X^*))$ and by Hilbert's Theorem 90, they are coboundaries, i.e.\ $b_w(\chi)= \beta(\chi)^{-1}\beta(^w\chi)$ for some $\beta\in \GL_n(\CC(\delta^{-\frac{1}{2}}A_X^*))$. Then $\beta(\chi)$ is the transition matrix between the basis $(S_{\tilde\chi})_{\tilde\chi}$ and the desired basis $(S^i_\chi)_i$. 

This basis $S^i_\chi$ has the problem that it is not explicit. If there are no orbits $\YY$ of maximal rank and simple roots $\alpha$ such that $(\YY,\alpha)$ is of type $N$ then we simply denote by $(S^i_\chi)_i$ the basis consisting of the morphisms $S_{\tilde\chi}$.

We can now state the main representation-theoretic result of this paper:

\begin{theorem}\label{mainthm}
 Assume that $\mathring X$ carries a $B$-invariant measure (\S\ref{invariantmeasure}) and let $\mathcal B_X$ denote the image of $\delta^{-\frac{1}{2}}A_X^*$ on the ``unramified'' Bernstein variety (\S\ref{ssBernstein}). By Theorem \ref{existence}, every irreducible $\pi\in\mathcal S_\ur$ admitting a non-zero quotient: $C_c^\infty(X)\to\pi$ must lie over $\mathcal B_X$. With the possible exception of a set of $\pi$'s (resp.\ $\chi$'s) lying over a proper closed subvariety of $\mathcal B_X$, the following are true:
\begin{enumerate}
 \item\label{mainthmfirst} Every $\pi$ is isomorphic to $\Ind_{P(X)}^G (\chi\delta^{\frac{1}{2}})$ for $\chi\in\delta^{-\frac{1}{2}}A_X^*$.
 \item\label{mainthmsecond} Every quotient $C_c^\infty(X)\to\pi$ is obtained as a linear combination of specializations of the morphisms $S^i_\chi$, $i=1,\dots,|H^1(k,A_X)|$.
 \item\label{mainthmthird} The quotients $S^i_\chi$, for fixed $\chi$, are linearly independent. The quotients $S^i_{\chi_1}$ and $S^j_{\chi_2}$ are isomorphic if and only if $\chi_1={^w\chi_2}$ for some $w\in W_X$ and $i=j$.
 \item\label{mainthmfourth} $\dim\Hom(C_c^\infty(X),\pi)=(\mathcal N_W(-\rho+\mathfrak a_X^*):W_X)\times{|H^1(k,\AA_X)|}$.
\end{enumerate}
\end{theorem}

\begin{proof}
\ref{mainthmfirst}) is Lemma \ref{regular} and Corollary \ref{irreducible}.

\ref{mainthmsecond}) follows from the fact that generically quotients into irreducible $\pi\in\mathcal S_\ur$ are obtained by (linear combinations of) specializations of morphisms in $\underline S_\chi^Y$ for $Y$ of maximal rank, Corollary \ref{cor1} and the fact that Knop's action is transitive on orbits of maximal rank.

\ref{mainthmthird}) is a consequence of Corollary \ref{cor1}, the definition of the $S_\chi^i$'s above.

\ref{mainthmfourth}) follows immediately from \ref{mainthmsecond}) and \ref{mainthmthird}).

\end{proof}

Let us now compare our results with a few well-known examples:

\begin{example}
The spherical variety $\XX=\UU\backslash\GG$. It is known that the little Weyl group of a horospherical variety is trivial (and vice versa: if the little Weyl group of a spherical variety is trivial then the variety is horospherical) and it is easy to check that $\Adm_X=A_X^*=A^*$. Therefore, our results translate to the fact that all irreducible representations in the unramified spectrum appear, at least generically, but the generic multiplicity is equal to the order of the Weyl group. (This is, of course, expressed by the isomorphisms $I(\chi)\simeq I(^w\chi)$ for generic $\chi$.)
\end{example}

\begin{example}
The subgroup $\HH=\GGL_n^{\rm{diag}}$ of $\GG=\GGL_n\times\GGL_n$. In this case $\Adm_X=A_X^*$ is equal to $A_{\GL_n}^*$ embedded in $A^*$ as $a\mapsto\it{diag}(a,a^{-1})$ and $W_X=W_{\GL_n}^{\rm{diag}}$. Therefore, generically in the unramified spectrum, $X=H\backslash G$ distinguishes (with multiplicity one) irreducible representations of the form $\tau\otimes\tilde \tau$, where $\tau$ is an irreducible representation of $\GL_n$ and $\tilde \tau$ denotes its contragradient. This, of course, holds not only generically and not only for the unramified spectrum.
\end{example}

\begin{example}
The space $\XX=\operatorname{\mathbf{Mat}}_n$ under the $\GG=\GGL_n\times\GGL_n$ action by multiplication on the left and right. The open $\GG$-orbit is equal to the spherical variety of the previous example, therefore the generic description of the unramified spectrum is identical to the previous case.
\end{example}

\begin{remark}
 As follows immediately from Theorem \ref{mainthm}, the generic multiplicity may be greater than 1 (i.e.\ the Gelfand condition may fail to hold) for two reasons: The $k$-points of the open $\BB$-orbit split into several $B$-orbits; or the little Weyl group of $\XX$ does not coincide with the normalizer of $-\rho+\mathfrak a_X^*$. In addition to the simple $\SL_2$-examples that we have seen, we mention another instance of the former:
\end{remark}

\begin{example}
In \cite{HiSp}, Y.\ Hironaka examines ${\Ssp}_4$ as a spherical homogeneous ${\Ssp}_4\times({\Ssp}_2)^2$-space over a local non-archimedean field. It is discovered that the generic multiplicity is equal to the order of $k^\times/(k^\times)^2$; this is due to the splitting of the $\mathring X$ in $B$-orbits. (Compare also the explicit computations in loc.\ cit.\ with our treatment of ``Case T'' in Proposition \ref{step1}.
\end{example}

The non-coincidence of $W_X$ with $\mathcal N(-\rho+\mathfrak a_X^*)$ is very common in parabolically induced examples, since, as we already mentioned, the little Weyl group of the parabolically induced spherical variety is equal to the little Weyl group of the original spherical variety for the Levi. The example of $\UU\backslash\GG$, mentioned above, is an instance of this. However, parabolically induced spherical varieties do not exhaust the list of such examples:

\begin{example}
The group $\SSL_2\times\SSL_2$ embeds naturally in $\GG=\Ssp_4$ as $\Ssp_2\times\Ssp_2$. Let $\HH$ be the $\Gm\times\SSL_2$ subgroup thereof (where $\Gm$ is a maximal split torus in $\SSL_2$). It is easy to see that $\Adm_X=A_X^*=A^*$, however it is known that $W_X$ is not the whole Weyl group, but a subgroup of $W$ of index 2. 
\end{example}

\subsection{Parabolic induction with an additive character}\label{ssparindthm}

In applications one often comes across representations induced from ``parabolically induced'' spherical subgroups, but not from the trivial (or the modulus of an algebraic) character of those subgroups but from a complex character of its unipotent radical. 

\begin{example}
The \emph{Whittaker model} is the line bundle $\mathcal L_\Psi$ over $U\backslash G$, where $\Psi: U\to\CC^\times$ is a \emph{generic character} of $U$; this means that, if we identify the abelianization of $\UU$ with the direct product of the one-parameter subgroups $\UU_\alpha$, for $\alpha$ ranging over all simple positive roots, then $\Psi=\psi\circ\Lambda$, where $\Lambda: \UU\to\Ga$ is a functional which does not vanish on any of the $\UU_\alpha$ and $\psi$ is a nontrivial complex character of $\Ga(k)=k$. (Equivalently, $\Lambda$ lies in the open $\AA$-orbit on $\Hom(\UU,\Ga)$.) Hence, the Whittaker model is parabolically induced from the trivial subgroup of $A$; if $\Psi$ were the trivial character, then its spectrum would only contain representations whose Jacquet module with respect to $U$ is non-trivial, and with generic multiplicity equal to the order of the Weyl group. On the contrary, for $\Psi$ a generic character, the spectrum is known to be much richer (e.g.\ it contains all supercuspidals), and \emph{multiplicity free} for every (not only generic) irreducible representation.
\end{example}

\begin{example}
The \emph{Shalika model} for $\GL_{2n}$ is the line bundle $\mathcal L_\Psi$ over $H\backslash G$, where $\HH$ is parabolically induced from the maximal parabolic $\PP$ with Levi $\LL=\GGL_n\times\GGL_n$ and the spherical subgroup $\MM=\GGL_n^{\rm{diag}}$ thereof; and $\Psi$ is the character  $\psi(\tr(X))$ of $U_P$, where $\psi$ is as above a complex character of $k$ and $X\in\it{Mat}_n(k)$ under the isomorphism $\UU_P\simeq \operatorname{\mathbf{Mat}}_n$. It is known that the Shalika model, too, is multiplicity-free, and it distinguishes lifts from $\operatorname{SO}_{2n+1}$.
\end{example}

We will see that even those cases can be linked to Knop's theory -- more precisely, to an extension of the Weyl group action to non-spherical varieties.

Let $\HH=\MM\ltimes\UU_\PP$ be a parabolically induced spherical subgroup of $\GG$, with notation as in \S \ref{ssparind}. 
Let $\Psi:U_P\to\CC^\times$ be a character. Any such character of $U_P$ factors through a morphism: $\Lambda:\UU_P\to \Ga$, composed with a complex character $\psi$ of $\Ga(k)=k$. Now, assume that $\Lambda$ is normalized by $\MM$, and by abuse of notation use the same letter to denote the induced morphism: $\HH\to\Ga$. Let $\HH_0=\ker\Lambda$; the variety $\HH_0\backslash\GG$ is the total space of a $\Ga$-torsor over $\HH\backslash\GG$ (no longer spherical), and the map $$\lambda:\HH_0\backslash\GG \to \HH\backslash \GG$$ is surjective on $k$-points for the usual reasons. One is interested in the space $C_c^\infty(\HH\backslash\GG,\mathcal L_\Psi)$, i.e.\ the space of smooth complex functions on $\HH_0\backslash\GG$ which satisfy $f(h\cdot x)=\Psi (h) f(x)$ for $h \in \HH/\HH_0 (k)$ and such that $\lambda(\supp(f))$ is compact. 

By repeating exactly the same Mackey-theoretic arguments that we used before, one sees directly that the $k$-rational $\BB$-orbits on $\HH\backslash\GG$ which give rise to a morphism into $I(\chi)$, for some unramified character $\chi$, are those represented by elements $\xi$ such that 
\begin{equation}\label{parabcond}
\HH\cap{^\xi\BB}\subset\HH_0.
\end{equation}
One sees also that if an orbit $\YY$ satisfies this condition, then one can define a rational family $S_\chi^Y$ of morphisms into $I(\chi)$ for exactly the same $\chi$'s as before. Also, by the description of $\BB$-orbits in \S \ref{ssparind}, one sees that the open $\BB$-orbit satisfies (\ref{parabcond}). Denote by $\mathfrak B_{00}^\Lambda$ the set of orbits of maximal rank which satisfy (\ref{parabcond}). Our goal is to describe the unramified quotients of $C_c^\infty(X,\Psi)$ in a similar manner as we did with $C_c^\infty(X)$; more precisely, we will link it with Knop's Weyl group action on $\HH_0\backslash\GG$.

Since the latter space is not spherical, we need to revisit Knop's theory and recall the necessary facts regarding its extension to non-spherical varieties.
The \emph{complexity} of a variety $\YY$ with a $\BB$-action is defined as $c(\YY)=\{\rm{max}_{y\in\YY}\rm{codim} (y\BB)\}$. Let $\XX$ be a $\GG$-variety, not necessarily spherical. We have $c(\XX)=0$ if and only if $\XX$ is spherical. We let $\mathfrak B_{0}(\XX)$ denote the set of closed, irreducible, $\BB$-stable subsets with complexity equal to the complexity of $\XX$. Then Knop defines an action of the Weyl group $W$ on $\mathfrak B_0$ -- it leaves stable the subset $\mathfrak B_{00}$ of those $\BB$-stable subsets whose general $\BB$-orbit is of maximal rank. In the case of spherical varieties, this action coincides with the one that we discussed, and $\mathfrak B_{00}$ is in bijection with the set of $\BB$-orbits of maximal rank, hence the use of the same symbol to denote those. To see how the action is defined in the general case, one repeats the same steps, by letting $\PP_\alpha$ act on the $\BB$-stable set $\YY\in\mathfrak B_{0}$, examining the image of a general stabilizer $\PP_y$ in $\Aut(\BB\backslash\PP_\alpha)\simeq \PPGL_2$ and considering cases. Additionaly to the cases that we saw in the spherical case, one now has the case $\FF\backslash\PPGL_2$, where $\FF$ is a finite subgroup. But in that case, there is only one closed, irreducible, $\BB$-stable subset of complexity equal to the complexity of $\FF\backslash\PPGL_2$ (namely, the space $\FF\backslash\PPGL_2$ itself), and the corresponding element of the Weyl group will by definition fix it.

Now let us return to our parabolically induced spherical variety. Let us consider inverse images of $\BB$-orbits under $\lambda$. The set $\{ \overline{\lambda^{-1}\YY} | \YY\in\BB_{00}(\HH\backslash\GG)\}$ is precisely the set of closed, irreducible, $\BB$-stable subsets of $\HH_0\backslash\GG$ whose generic $\BB$-orbit has maximal rank. Which of those belong to $\mathfrak B_{00}(\HH_0\backslash\GG)$? One sees easily that, for $\YY$ a $\BB$-orbit on $\HH\backslash\GG$,  $\overline{\lambda^{-1}\YY}$ has complexity 1 (the complexity of $\HH_0\backslash\GG$) if and only if $\lambda^{-1}\YY$ is not a single $\BB$-orbit, which is the case if and only if (\ref{parabcond}) is satisfied. Therefore \emph{we have a natural isomorphism of sets}: $\mathfrak B_{00}^\Lambda (\HH\backslash\GG) \simeq \mathfrak B_{00}(\HH_0\backslash\GG)$. The Weyl group action on the latter induces a Weyl group action on the former, which differs from the action of $W$ on $\mathfrak B_{00}(\XX)$. Then:

\begin{theorem}
In the above setting, let $\underline S_\chi^Y$ denote the family of morphisms into $I(\chi)$ defined by the $k$-rational orbit $\YY\in\mathfrak B_{00}^\Lambda(\XX)$  and let ${^w\YY}$ denote the image of $\YY$ under the Weyl group action on $\mathfrak \BB_{00}^\Lambda(\XX)$.

Then $T_w\circ \underline S_\chi^{\mathring X} \ne 0$ if and only if $w\in [W/W_{P(X)}]$. For $w\in [W/W_{P(X)}]$, $T_w\circ \underline S_\chi^{\mathring X} = \underline S_{^w\chi}^{{^w\mathring X}}$. If there are no orbits $\YY$ of maximal rank and simple roots $\alpha$ such that $(\YY,\alpha)$ is of type $N$, then  $T_w\circ  S_{\tilde\chi}^{\mathring X} \sim \underline S_{^w\tilde\chi}^{{^w\mathring X}}$ for every $w\in [W/W_{P(X)}]$.
\end{theorem}

Theorem \ref{mainthm} extends verbatim to this setting, with $W_X$ the ``little Weyl group'' of $\HH_0\backslash\GG$.
The proof is similar to the case of $\Psi=$ trivial and is omitted.


\section{Unramified vectors and endomorphisms}\label{secEnd}

\subsection{Spectral support}

Since the results of this paper are all stated for ``generic'' quotients of the ``unramified'' Bernstein component, it is natural to ask to what extent those ``generic'' quotients are enough to characterize a vector in our representation. Given a smooth representation $\pi$, let us call \emph{spectral support} (or simply \emph{support}) of $\pi$ its support as a module for the Bernstein center $\mathfrak z(\mathcal S)$. 
In other words, it is the subvariety of the Bernstein variety corresponding to the ideal of $\mathfrak z(\mathcal S)$ which annihilates $\pi$. Given a vector $v\in \pi$, we will call \emph{(spectral) support} of $v$ the support of the smallest subrepresentation of $\pi$ containing $v$. Our question can be reformulated as follows: To what extent is the spectral support of a vector $v\in C_c^\infty(X)_\ur$ (or one of its quotients) equal to the image $\mathcal B_X$ of $\delta^{-\frac{1}{2}}A_X^*$ in $A^*/W$ (the unramified component of the Bernstein variety)? We will say that $v$ is of ``generic support'' if its support is equal to $\mathcal B_X$.

It is easy to see that not all vectors in $C_c^\infty(X)_\ur$ have generic support in general. For instance, let $\XX=\TT\backslash\PPGL_2$ as in Example \ref{ex2}. Recall our description of its Jacquet module: We can have $\phi\in C_c^\infty(X)$ whose image in the Jacquet module is non-zero, but is zero when restricted to $k^\times$. In fact, we can generate $\phi$ as follows: Choose a suitable $\phi_1$ supported in a neighborhood of the divisor $Y_1$ (in the notation of \ref{ex2}), and apply the automorphism ``$w$'' to it (action of the non-trivial Weyl group element on the left). Let $\phi=\phi_1-{^w\phi_1}$. Since ``$w$'' is $G$-equivariant, we will have $R(g) \phi = R(g) \phi_1 - {^w(R(g)\phi_1)}$, and therefore the image of all translates $\phi$ in the Jacquet module will be supported at the ``double origin'' . As a result, the support of $\phi$ is not generic.

Let $K$ be a hyperspecial maximal compact subgroup of $G$. The following theorem gives an assertive answer to our question for $K$-invariant vectors.
For what follows we will be denoting by $S_{\tilde\chi}, S_{\tilde\chi,\zeta}$ the $B$-equivariant functionals into $\CC_{\chi\delta^{\frac{1}{2}}}$ defined by the open $\BB$-orbit $\mathring\XX$.

The following theorem is clearly false in the case of anomalous non-homogeneous varieties such as that of Example \ref{ugly}. Therefore, \emph{for the rest of the paper we re-define $X$ to mean the Hausdorff closure of $\mathring X$} (cf.\ Lemma \ref{smorbitlemma}):

\begin{theorem}\label{support}
The support of every $\phi\in C_c^\infty(X)^K$ is generic. In fact, if $S_{\tilde\chi}(\phi) = 0$ (as a functional) for almost every $\tilde\chi$ then $f=0$.
\end{theorem}

\begin{proof}
The second statement, although it appears stronger, is actually equivalent. First, by our main theorem if $S_{\tilde\chi}(\phi)=0$ as a \emph{morphism} into $I(\chi)$ (\emph{not} as a functional) then $S_{\tilde\chi}^Y(\phi)=0$ for every orbit $\YY$ of maximal rank. Moreover $I(\chi)$ contains a unique (up to scaling) non-zero $K$-invariant vector whose value at 1 is non-zero, so to say that the value of the functional $S_{\tilde\chi}$ applied to $\phi$ is zero, for a $K$-invariant $\phi$, is the same as saying that $S_{\tilde\chi}(\phi)=0$ as an element of $I(\chi)$.

Assume $S_{\tilde\chi}(\phi)=0$ for almost every $\tilde\chi$. By Lemma \ref{enough}, it suffices to show that the functionals $S_{\tilde\chi}^Y$, for $\YY$ of smaller rank, vanish on $\phi$ and all its translates. We will consider two cases separately: $\YY$ a Borel-orbit in the open $\GG$-orbit, and $\YY$ a Borel-orbit in a different $\GG$-orbit.

\subsubsection{Proof within the open $\GG$-orbit.} We use induction on the dimension of $\YY$ to show that if $S_{\tilde\chi}^Z(\phi)$ vanishes identically for all orbits $\ZZ$ of dimension larger than the dimension of $\YY$ then $S_{\tilde\chi}^Y(\phi)$ vanishes, too.

First, if  $\YY$ is raised of type U by some simple root $\alpha$ to an orbit $\ZZ$ then we can apply Theorem \ref{maintool} to deduce the vanishing of $S_{\tilde\chi}^Y(\phi)$ from the vanishing of $S_{\tilde\chi}^Z(\phi)$.

Now, assume that there exists a simple root $\alpha$ raising $\YY$ of type T or N. Again by Theorem \ref{maintool}, we would like to deduce the vanishing of $S_{\tilde\chi}^Y(\phi)$ from the vanishing of $S_{\tilde\chi}^{^{w_\alpha}Y}(\phi)$ (resp.\ $S_{\tilde\chi,\zeta}^{^{w_\alpha}Y}(\phi)$ in case N). The point to be careful about here is that $\Adm_Y$ may be contained in one of the divisors where $T_{w_\alpha}$ is not an isomorphism. However, by the remark in \S\ref{sssadmissible}, that divisor can only be $\mathcal Q_{-\alpha}$ and that implies that $T_{w_\alpha}$ does not ``kill'' $I(\chi)^K$.

This finishes the case where $\YY$ is contained in the open $\GG$-orbit.

\subsubsection{Proof on smaller $\GG$-orbits.} Let $\YY$ be a $k$-rational Borel orbit, belonging to a non-open $\GG$-orbit $\ZZ$. Since $\ZZ$ itself is a spherical variety, it suffices by the proof of the previous case to assume that $\YY=\mathring\ZZ$ (the open $\BB$-orbit in $\ZZ$). Also, we may inductively assume that the theorem has been proven for all larger $\GG$-orbits containing $\ZZ$ in their closure. Let $Z_1$ be a $B$-orbit on $\mathring Z$ and let $V_1$ be a $B$-orbit on a larger $G$-orbit as in Proposition \ref{pickup}, containing $Z_1$ in its closure. (We are using here the fact that $X$ was redefined as the Hausdorff closure of $\mathring X$.) By Proposition \ref{pickup} we know that $S_\chi^{Z_1}$ is a residue of $S_\chi^{V_1}$, therefore since the latter is identically zero on $\phi$ the former must be, also. This completes the proof of the theorem.

\end{proof}

\subsection{The Hecke module of unramified vectors} The result of the previous section allows us to present a weak analog of the main result of \cite{GN1,GN2}, namely a description of the Hecke module of $K$-invariant vectors. (These vectors are commonly called ``spherical'', but to avoid a double use of this word we will only call them ``unramified''.) Notice that $C_c^\infty(X)^K\subset C_c^\infty(X)_\ur$. 

Let $\mathcal H(G,K)$ denote the convolution algebra of $K$-biinvariant measures on $G$. Recall the Satake isomorphism: $\mathcal H(G,K)\simeq \CC [A^*]^W$. Since all vectors in $C_c^\infty(X)_\ur$ have spectral support over the image of $\delta^{-\frac{1}{2}}A_X^*$ in $A^*/W$, $\mathcal H(G,K)$ acts on $C_c^\infty(X)^K$ through the corresponding quotient, which will be denoted by $\mathcal H_X$. Let $\mathcal K_X$ denote the fraction field of $\mathcal H_X$, hence naturally: $\mathcal K_X\simeq \CC(\delta^{-\frac{1}{2}}A_X^*)^{\mathcal N_W(-\rho+\mathfrak a_X^*)}$.

\begin{theorem}\label{Ktheorem}
The space $C_c^\infty(X)^K$ is a finitely-generated, torsion-free module for $\mathcal H_X$. 

Moreover, we have: $C_c^\infty(X)^K\otimes_{\mathcal H_X} \mathcal K_X \simeq \left(\CC(\delta^{-\frac{1}{2}}A_X^*)^{W_X}\right)^{|H^1(k,\AA_X)|}$.
\end{theorem}

\begin{remark}
The isomorphism above is not a canonical one, since it depends, as we shall see, on the choice of one $K$-invariant vector.
\end{remark}

\begin{proof}
The fact that $C_c^\infty(X)^K$ is torsion-free over $\mathcal H_X$ follows from Theorem \ref{support}.

Now let $(S_\chi^i)_i$ denote the operators of Theorem \ref{mainthm} (in particular, they are equal to the operators $S_{\tilde\chi}^{\mathring X}$ if  there are no pairs $(\YY,\alpha)$ where $\YY$ is a $\BB$-orbit of maximal rank, $\alpha$ is a simple root and $(\YY,\alpha)$ is of type N). 
Fix a $\phi_0\in C_c^\infty(X)^K$ and define a map $C_c^\infty(X)^K\to \left(\CC(\delta^{-\frac{1}{2}}A_X^*)\right)^{|H^1(k,\AA_X)|}$ by 
$$\phi\mapsto \left (\frac{S_{\chi}^i(\phi)}{S_{\chi}^i(\phi_0)}\right )_{i}.$$

We claim, first, that the image of this map lies in the $W_X$-invariants. Indeed, by Theorem \ref{mainthm}, for $w\in W_X$ the quotients $S_\chi^i$ and $S_{^w\chi}^i$ are isomorphic (through $T_w$). Since there is a unique line of $K$-invariant vectors in $I(\chi)$, it follows that $S_\chi^i(\phi)=c(\chi) \cdot S_\chi^i(\phi_0)$ for some constant $c(\chi)$, rational in $\chi$, and since $S_\chi^i$ and $S_{^w\chi}^i$ are isomorphic, it follows that that this constant is the same for $\chi$ and for $^w\chi$. This proves that we have a map: $C_c^\infty(X)^K \to \left(\CC(\delta^{-\frac{1}{2}}A_X^*)^{W_X}\right)^{|H^1(k,\AA_X)|}$. 

We have shown in Theorem \ref{support} that for a nonzero $\phi\in C_c^\infty(X)^K$ it is not possible that $S_{\chi}^i(\phi)=0$ for every $\chi,i$ -- this establishes that the map is injective.

We prove surjectivity of the map when tensored with $\mathcal K_X$. Notice that the space of morphisms: $C_c^\infty(X)\to I(\chi)$ has generically dimension equal to $r:=(\mathcal N(\rho+\mathfrak a_X^*):W_X)\cdot |H^1(k,\AA_X)|$. Suppose that $C_c^\infty(X)^K\otimes_{\mathcal H_X} \mathcal K_X$ had smaller dimension over $\mathcal K_X$. Then the basis $(S_\chi^i)_i$ of the space of morphisms: $C_c^\infty(X)\to I(\chi)$ satisfies a linear relation: $\sum_i c_i(\chi) S_\chi^i = 0$ (with $c_i(\chi)$ rational)  when restricted to the subspace generated by $C_c^\infty(X)^K$. Given $\phi\in C_c^\infty(X)_\ur$, I claim that $\sum_i c_i(\chi) S_\chi^i (\phi) = 0$ for generic (and hence every) $\chi$. Indeed, for generic $\chi$ the image of the $S_\chi^i$'s in $I(\chi)$ is irreducible and unramified, hence if $\sum_i c_i(\chi) S_\chi^i (\phi)$ is non-zero then one of its $G$-translates, when convolved with the characteristic function of $K$, should be non-zero. Since the $S_\chi^i$'s are $G$-equivariant and $R(Kg)\phi \in C_c^\infty(X)^K$, we have: $R(Kg)\sum_i c_i(\chi) S_\chi^i (\phi) = \sum_i c_i(\chi) S_\chi^i (R(Kg)\phi) = 0$ by assumption. It follows that $\sum_i c_i(\chi) S_\chi^i =0$ on the whole space, contradicting what we know about the dimension of the space of morphisms into $I(\chi)$. This proves the stated isomorphism.

Finally, recall from \S \ref{sspoles} that we may regularize the $S_\chi^i$'s (by multiplying by a suitable regular function of $\chi$) so that they are regular for every $\chi$. Then $\phi\mapsto (S_\chi^i(\phi))_i$ defines an injection $C_c^\infty(X)^K\hookrightarrow \CC[\delta^{-\frac{1}{2}}A_X^*]^{|H^1(k,\AA_X)|}$. Since the latter is a finitely generated $\mathcal H_X$-module, it follows that $C_c^\infty(X)^K$ is also finitely generated.

\end{proof}

\subsection{A commutative ring of endomorphisms}

In this section we assume, for simplicity, that $H^1(k,\AA_X)=\{1\}$, i.e.\ each $\BB$-orbit of maximal rank contains only one rational $B$-orbit.

\subsubsection{Definition of the map.}
Let $D\in \End_{\mathcal H(G,K)}(C_c^\infty(X)^K)$. It induces an endomorphism of $C_c^\infty(X)^K\otimes_{\mathcal H_X} \mathcal K_X \simeq \left(\CC(\delta^{-\frac{1}{2}}A_X^*)^{W_X}\right)$ which is $\mathcal K_X$-linear. If it is also $\CC(\delta^{-\frac{1}{2}}A_X^*)^{W_X}$-linear (in other words, if $S_\chi\circ D\sim S_\chi$ on $C_c^\infty(X)^K$), then we will call $D$ ``geometric''. Of course, if $(C_c^\infty(X))_\ur$ is generically multiplicity-free (i.e.\ $W_X = \mathcal N_W(-\rho+\mathfrak a_X^*)$), then every endomorphism is geometric, but this will not be the case in general. The map $D\mapsto c(\chi)$, where $c(\chi)$ is given by the relation $S_\chi\circ D = c(\chi) S_\chi$, defines a ring homomorphism $\End_{\mathcal H(G,K)}(C_c^\infty(X)^K)^{\geom} \to \CC(\delta^{-\frac{1}{2}}A_X^*)$. In fact, by Theorem \ref{maintool}, the image lies in invariants of the little Weyl group $W_X$. Moreover, since by Theorem \ref{Ktheorem} $C_c^\infty(X)^K$ is a finitely generated $\mathcal H_X$-module, every $\mathcal H_X$-algebra of endomorphisms is a finitely generated module over $\mathcal H_X$; and since the integral closure of $\mathcal H_X$ in $\CC(\delta^{-\frac{1}{2}}A_X^*)^{W_X}$ is $\CC[\delta^{-\frac{1}{2}}A_X^*]^{W_X}$ (the variety $A_X^*/W_X$ is normal), it follows that the image of the above homomorphism must lie in $\CC[\delta^{-\frac{1}{2}}A_X^*]^{W_X}$. We conjecture that the image is the whole ring:

\begin{conjecture}\label{conjecture}
There is a canonical isomorphism: $$\left(\End_{\mathcal H(G,K)}C_c^\infty(X)^K\right)^{\geom} \simeq \CC[\delta^{-\frac{1}{2}}A_X^*]^{W_X}.$$
\end{conjecture}

The reason that we believe that these endomorphisms exist in general is the following analogy with invariant differential operators on $\XX$:

As was proven by F.~Knop in \cite{KnHC}, the algebra of invariant differential operators on a spherical variety (over an algebraically closed field $K$ in characteristic zero) is commutative, and isomorphic to $K[\rho+\mathfrak a_X^*]^{W_X}$. Here, $\mathfrak a_X^*$ is isomorphic to the Lie algebra of what we denote by $A_X^*$. This generalizes the Harish-Chandra homomorphism (if we regard the group $\GG$ as a spherical $\GG\times \GG$ variety), and in fact the following diagram is commutative:

\begin{equation} \begin{CD}
\mathfrak{z}(\GG) @>>> \mathfrak{D}(\XX)^\GG \\
@|                       @| \\
K[\mathfrak a^*]^W @>>> K[\rho+\mathfrak a_X^*]^{W_X}
\end{CD}
\end{equation}

What we propose is a similar diagram for the $p$-adic group, which on the left side will have the Satake isomorphism for $\mathcal H(G,K)$ (or equivalently, the unramified factor of the Bernstein centre) and on the right side the ``geometric endomorphisms'' that we defined above, which should be viewed as an analog for the invariant differential operators:

\begin{equation} \begin{CD}
\mathcal H(G,K) @>>> \left(\End_{\mathcal H(G,K)}C_c^\infty(X)^K\right)^{\geom} \\
@|                @|  \\
\CC[A^*]^W @>>> \CC[\delta^{-\frac{1}{2}}A_X^*]^{W_X}
\end{CD}
\end{equation}

\begin{remark}
The reader should not be confused by the fact that in Knop's result the Lie algebra $a_X^*$ is offset by $\rho$ while in ours the torus $A_X^*$ is offset by $\delta^{-\frac{1}{2}}$, which is equal to $-\rho$: The discrepancy is a matter of definitions, and to fix it we could have denoted by $I(\chi^{-1})$ what we denoted by $I(\chi)$ -- but of course this would contradict the conventions in the literature.
\end{remark}
Our interest in this conjecture comes from the fact that the analog of invariant differential operators suggests the possibility of a ``geometric'' construction of these endomorphisms, while spectral methods do not seem to suffice in general.

However, it is easy to prove the conjecture in the cases that $X$ is generically multiplicity-free or parabolically induced from a multiplicity-free one. Then these endomorphisms will be provided by the Hecke algebra of $G$ or, respectively, of a Levi subgroup acting ``on the left'':

\begin{theorem}\label{conjpf}
Conjecture \ref{conjecture} is true if: 
\begin{enumerate}
\item
the unramified spectrum of $X$ is generically multiplicity-free, in which case the geometric endomorphisms are all the endomorphisms of $C_c^\infty(X)^K$, or
\item
the spherical variety $X$ is ``parabolically induced'' from a spherical variety whose unramified spectrum is generically multiplicity-free.
\end{enumerate}
\end{theorem}

\begin{proof}
In the first case, $\CC[A^*]^W$ surjects onto $\CC[\delta^{-\frac{1}{2}}A_X^*]^{W_X}$. The claim that these are all the endomorphisms follows from Theorem \ref{support}.

In the second case, as we saw in \S \ref{ssparind}, the subtorus $A_X^*$ and the Weyl group $W_X$ coincide with those associated to the corresponding spherical variety of the Levi. Hence $\CC[\delta^{-\frac{1}{2}}A_X^*]^{W_X}$ is surjected upon by the Bernstein centre of $L$ acting ``on the left''.
\end{proof}



\begin{thebibliography}{KKV89}

\bibitem[At70]{At} M.\ F.\ Atiyah, \emph{Resolution of singularities and division of distributions.} Comm.\ Pure Appl.\ Math.\ 23 1970 145--150. 

\bibitem[Be84]{BeCtr} J.\ N.\ Bernstein, \emph{Le ``centre'' de Bernstein.} Edited by P.~Deligne. Travaux en Cours,  Representations of reductive groups over a local field,  1--32, Hermann, Paris, 1984.

\bibitem[BG69]{BG} J.\ N.\ Bernstein and S.\ I.\ Gel'fand, \emph{Meromorphy of the function $P\sp{\lambda }$.} Funkcional.\ Anal.\ i Prilo\v zen.\  3  1969 no.\ 1, 84--85.

\bibitem[BZ76]{BZ} J.\ N.\ Bernstein and A.\ V.\  Zelevinsky, \emph{Representations of the group $\GL(n,F)$ where $F$ is a local non-Archimedean field.} Uspehi Mat.\ Nauk 31 (1976), no.\ 3(189), 5--70. 

\bibitem[Bo79]{BoAut} A.\ Borel, \emph{Automorphic $L$-functions.}  Automorphic forms, representations and $L$-functions (Proc.\ Sympos.\ Pure Math., Oregon State Univ., Corvallis, Ore., 1977), Part 2,  pp.\ 27--61, Proc. Sympos. Pure Math., XXXIII, Amer. Math. Soc., Providence, R.I., 1979.

\bibitem[Bo91]{Bo} A.\ Borel, \emph{Linear algebraic groups.} Second edition, Graduate Texts in Mathematics, 126, Springer-Verlag, New York, 1991.

\bibitem[BK99]{BK} A.\ Braverman and D.\ Kazhdan, \emph{On the Schwartz space of the basic affine space.}  Selecta Math.\ (N.S.)  5  (1999),  no.\ 1, 1--28.

\bibitem[Br86]{BrPr}  M.\ Brion, \emph{Quelques propri\'et\'es des espaces homog\`enes sph\'eriques.} Manuscripta Math.\  55  (1986),  no.\ 2, 191--198.

\bibitem[Br90]{BrGe} M.\ Brion, \emph{Vers une g\'en\'eralisation des espaces symm\'etriques.}  J.\ Algebra  134  (1990),  no.\ 1, 115--143.

\bibitem[Br01]{BrOr} M.\ Brion, \emph{On orbit closures of spherical subgroups in flag varieties.}  Comment.\ Math.\ Helv.\  76  (2001),  no.\ 2, 263--299.

\bibitem[Ca80]{Csph} W.\ Casselman, \emph{The unramified principal series of $p$-adic groups. I. The spherical function.}  Compositio Math.\ 40  (1980), no.\ 3, 387--406.

\bibitem[Cas]{Cas} W.\ Casselman, \emph{Introduction to the theory of admissible representations of $p$-adic reductive groups.} Draft, 1 May 1995, available at: \\ \texttt{http://www.math.ubc.ca/$\sim$cass/research/p-adic-book.dvi}.

\bibitem[Den85]{Den} J.\ Denef, \emph{On the evaluation of certain $p$-adic integrals.} S\'eminaire de th\'eorie des nombres, Paris 1983--84, 25--47, Progr.\ Math., 59, Birkhäuser Boston, Boston, MA, 1985. 

\bibitem[Des96]{Dh} B.\ Deshommes, \emph{Crit\`eres de rationalit\'e et application \`a la s\'erie g\'en\'eratrice d'un syst\`eme d'\'equations \`a coefficients dans un corps local.}  J.\ Number Theory  22  (1986),  no.\ 1, 75--114.\ 

\bibitem[GN1]{GN1} D.\ Gaitsgory and D.\ Nadler, \emph{Hecke operators on quasimaps into horospherical varieties.} Preprint, 2004, math.AG/0411266.

\bibitem[GN2]{GN2} D.\ Gaitsgory and D.\ Nadler, \emph{Spherical varieties and Langlands duality.} Preprint, 2006, math.RT/0611323.

\bibitem[Ga99]{Ga} P.\ Garrett, \emph{Euler factorization of global integrals.}  Automorphic forms, automorphic representations, and arithmetic (Fort Worth, TX, 1996),  35--101, Proc.\ Sympos.\ Pure Math., 66, Part 2, Amer.\ Math.\ Soc., Providence, RI, 1999.

\bibitem[GPR87]{GPR} S.\ Gelbart, I.\ Piatetski-Shapiro and S.\ Rallis, \emph{Explicit constructions of automorphic $L$-functions.} Lecture Notes in Mathematics, 1254, Springer-Verlag, Berlin, 1987.

\bibitem[HW93]{HW} A.\ G.\ Helminck and S.\ P.\ Wang, \emph{On rationality properties of involutions of reductive groups.}  Adv.\ Math.\  99  (1993),  no.\ 1, 26--96.

\bibitem[HiH64]{Hi}  H.\ Hironaka, \emph{Resolution of singularities of an algebraic variety over a field of characteristic zero.\ I, II.}  Ann.\ of Math.\ (2) 79 (1964), 109--203; ibid.\ (2)  79  1964 205--326.\ 

\bibitem[HiY99]{HiHe} Y.\ Hironaka, \emph{Spherical functions and local densities on Hermitian forms.}  J.\ Math.\ Soc.\ Japan  51  (1999),  no.\ 3, 553--581.

\bibitem[HiY05]{HiSp} Y.\ Hironaka, \emph{Spherical functions on ${\rm Sp}\sb 2$ as a spherical homogeneous ${\rm Sp}\sb 2\times({\rm Sp}\sb 1)\sp 2$-space.}  J.\ Number Theory  112  (2005),  no.\ 2, 238--286.

\bibitem[Ig00]{Ig} J-i.\ Igusa, \emph{An introduction to the theory of local zeta functions.} AMS/IP Studies in Advanced Mathematics, 14.\ American Mathematical Society, Providence, RI; International Press, Cambridge, MA, 2000.

\bibitem[Ja97]{Ja} H.\ Jacquet, \emph{Automorphic spectrum of symmetric spaces.} Representation theory and automorphic forms (Edinburgh, 1996),  443--455, Proc.\ Sympos.\ Pure Math., 61, Amer.\ Math.\ Soc., Providence, RI, 1997.

\bibitem[Kn90]{KnWe}  F.\ Knop, \emph{Weylgruppe und Momentabbildung.}  Invent.\ Math.\  99  (1990),  no.\ 1, 1--23.

\bibitem[Kn91]{KnLV}  F.\ Knop, \emph{The Luna-Vust theory of spherical embeddings.}  Proceedings of the Hyderabad Conference on Algebraic Groups (Hyderabad, 1989),  225--249, Manoj Prakashan, Madras, 1991.\ 

\bibitem[Kn94a]{KnAs} F.\ Knop, \emph{The asymptotic behavior of invariant collective motion.}  Invent.\ Math.\  116  (1994),  no.\ 1-3, 309--328.

\bibitem[Kn94b]{KnHC} F.\ Knop, \emph{A Harish-Chandra homomorphism for reductive group actions.}  Ann.\ of Math.\ (2)  140  (1994),  no.\ 2, 253--288.

\bibitem[Kn95a]{KnOrbits}  F.\ Knop, \emph{On the set of orbits for a Borel subgroup.}  Comment.\ Math.\ Helv.\  70  (1995),  no.\ 2, 285--309.

\bibitem[Kn95b]{KnR1} F.\ Knop, \emph{Homogeneous varieties for semisimple groups of rank one.}  Compositio Math.\  98  (1995),  no.\ 1, 77--89.

\bibitem[Kn96]{KnAu} F.\ Knop, \emph{Automorphisms, root systems, and compactifications of homogeneous varieties.} J.\ Amer.\ Math.\ Soc.\ 9 (1996), no.\ 1, 153--174. 

\bibitem[KKV89]{KKV} F.\ Knop, H.\ Kraft and T.\ Vust, \emph{The Picard group of a $G$-variety.}  Algebraische Transformationsgruppen und Invariantentheorie,  77--87, DMV Sem., 13, Birkh\"auser, Basel, 1989.

\bibitem[La]{La} E.\ Lapid, \emph{The relative trace formula and its applications.} Preprint, RIMS, Kyoto, 2005.

\bibitem[Of04]{Of}
 O.\ Offen, \emph{Relative spherical functions on $p$-adic symmetric spaces (three cases).}  Pacific J.\ Math.\  215  (2004),  no.\ 1, 97--149.\ 

\bibitem[PS75]{PS} I.\ Piatetski-Shapiro, \emph{Euler subgroups.}  Lie groups and their representations (Proc.\ Summer School, Bolyai J\'anos Math.\ Soc., Budapest, 1971),  pp.\ 597--620.\ Halsted, New York, 1975.\ 

\bibitem[Po86]{Po} V.\ L.\ Popov, \emph{Contractions of actions of reductive algebraic groups.} Mat.\ Sb.\ (N.S.) 130 (172) (1986), no.\ 3, 310--334, 431. 

\bibitem[Vi86]{Vi}\`E.\ B.\ Vinberg, \emph{Complexity of actions of reductive groups.} Funktsional.\ Anal.\ i Prilozhen.\  20  (1986),  no.\ 1, 1--13, 96.

\bibitem[VK78]{VK}\`E.\ B.\ Vinberg and B.\ N.\ Kimel'feld, \emph{Homogeneous domains on flag manifolds and spherical subsets of semisimple Lie groups.} Funktsional.\ Anal.\ i Prilozhen.\ 12  (1978), no.\ 3, 12--19, 96.

\bibitem[We82]{We} A.\ Weil, \emph{Adeles and algebraic groups.} Progr.\ Math.\ 23, Birkh\"auser, Boston, 1982.

\end{thebibliography}
\end{document}